# Mixed displacement-pressure-phase field framework for finite strain fracture of nearly incompressible hyperelastic materials

Fucheng Tian[a], Jun Zeng[a], Mengnan Zhang[a], Liangbin Li[a*]

a *National Synchrotron Radiation Laboratory, Anhui Provincial Engineering Laboratory of Advanced Functional Polymer Film, CAS Key Laboratory of Soft Matter Chemistry, University of Science and Technology of China, Hefei, 230026, China*

## Abstract

The favored phase field method (PFM) has encountered challenges in the finite strain fracture modeling of nearly or truly incompressible hyperelastic materials. We identified that the underlying cause lies in the innate contradiction between incompressibility and smeared crack opening. Drawing on the stiffness-degradation idea in PFM, we resolved this contradiction through loosening incompressible constraint of the damaged phase without affecting the incompressibility of intact material. By modifying the perturbed Lagrangian approach, we derived a novel mixed formulation. In numerical aspects, the finite element discretization uses the classical Q1/P0 and high-order P2/P1 schemes, respectively. To ease the mesh distortion at large strains, an adaptive mesh deletion technology is also developed. The validity and robustness of the proposed mixed framework are corroborated by four representative numerical examples. By comparing the performance of Q1/P0 and P2/P1, we conclude that the Q1/P0 formulation is a better choice for finite strain fracture in nearly incompressible cases. Moreover, the numerical examples also show that the combination of the proposed framework and methodology has vast potential in simulating complex peeling and tearing problems.

***Keywords***: Phase field model; Incompressible; Finite strain fracture; Hyperelastic materials

## 1. Introduction

With the traits of bearing large deformations and high recoverability, hyperelastic materials occupy a unique status in industrial and biomedical applications like seals, tires, artificial soft tissues, etc. Meanwhile, they also serve as an ideal vehicle for fundamental research such as material physics and nonlinear mechanics [1]. The dual demands of practicality and academic pursuit have boosted in-depth explorations on the service behavior of hyperelastic materials. Among them, finite strain, or large strain, fracture as the main failure channel of such materials, has intrigued widespread experimental and numerical research interest [2-6]. Yet the finite strain coupled discontinuous fracture, especially the tracking issue of the complex crack surfaces involved in three-dimensional fracture, places a tough obstacle for numerical treatments. To overcome this numerical challenge, the state-of-the-art diffusive crack models came into being, in which the sharp crack topology is approximated by the continuous scalar damage field [7-9]. As a prominent member of such models, the phase field model (PFM) has been sought after by the computational mechanics community in the past two decades [10, 11], due to its inherent features that can be naturally incorporated into the framework of continuum mechanics.

* Correspondence author: lbli@ustc.edu.cn (Liangbin Li)

The PFM prevailing in the computational mechanics community is rooted in the classic Griffith's theory [12], and shares some thoughts with those originated from Ginzburg-Landau phase transition theory [13-15]. Limiting our horizons to the former, the pioneering works are the variational approach to fracture established by Francfort and Marigo and the Ambrosio-Tortorelli regularized version reformulated by Bourdin et al. [16, 17]. Thenceforth, the research on the phase field methodology and its application has entered a flourishing stage [18-22]. Although the generic phase field model was developed for quasi-static brittle fracture, its application scenarios have been expanded to brittle cohesive fracture [23], dynamic fracture [24, 25], and finite-strain fracture [26-29], covering almost all material categories. Instead of presenting an exhaustive research list, we tend to review the studies on the fracture of hyperelastic materials. In this respect, Miehe et al. first presented a rate-independent phase field fracture (PFF) model for rubber-like polymer [30]. Loew et al. further extended the model to a rate-dependent type, which was validated by experiments on rubber fractures [27]. To better distinguish the contribution of tension and compression to finite-strain fracture, Tang et al. proposed a novel energy decomposition scheme for Neo-Hookean materials [31]. Moreover, Tian et al. and Peng et al. reformatted the PFM in the framework of edge-based smooth finite element method (ES-FEM), relieving the mesh distortion of large-strain fracture [32, 33]. Other associated researches involve hydrogel systems [34, 35], anisotropy [36], and composite materials [37]. Note, however, that most existing studies on the phase field modeling of hyperelastic materials focused on compressible cases, and rarely concerned incompressibility.

Although compressibility is commonly presumed in PFF modeling, we remark that numerous hyperelastic materials like rubber-like polymers and biological tissues can stand large strains without distinct volume changes. Therefore, these materials can be rationally treated as nearly or fully incompressible. In such circumstances, standard single-field displacement solutions suffer from the well-known volumetric locking problem [38, 39]. To alleviate or avert locking at finite deformation, many sophisticated numerical strategies and techniques have been developed, e.g., $\boldsymbol{F}$-bar and $J$-bar methods [40-42], enhanced-strain approach [43], mixed two-field displacement-pressure ($\boldsymbol{u}$-$p$) formulation and three-field displacement-pressure-Jacobian ($\boldsymbol{u}$-$p$-$J$) formulation [44, 45]. Regarding the specific pros and cons of these methodologies, the discussion is beyond the scope of this article and will not be further extended.

As stated earlier, incompressibility is an inescapable challenge in the modeling of most hyperelastic materials. However, the previous studies on the PFF modeling seldom touched on this issue, regardless of the small or the finite deformation. To our knowledge, Ma et al. first reported that the mixed $\boldsymbol{u}$-$p$ formulation was utilized to dispose of the phase field modeling of incompressible hyperelastic materials [46]. They modified the Lagrangian multiplier by introducing slightly compressible, thereby realizing simultaneous damage of bulk and shear modulus. Along this route, Li et al. pointed that the bulk modulus should be degraded faster than the shear modulus to free the damaged materials from the incompressible constraint [47]. However, this notion has not gained enough attention, owing to the short of solid theoretical foundation and detailed numerical implementation. More recently, Ye et al. proposed a new phase field model for nearly incompressible hyperelastic materials based on an enhanced three-field variational formulation [48]. They stated that introducing the assumed strain method also mediated the contradiction between energy decomposition and incompressible constraints. Despite the three-field formulation being mathematically elegant, we must recognize that the computational cost is also expensive, especially for PFF modeling.

Against this background, we endeavor to formulate a novel multi-field framework for modeling the fracture of nearly incompressible hyperelastic materials in finite strain scenarios. The underlying contradiction between incompressibility and the diffuse crack opening is elaborated from the perspective of physical crack topology. To resolve this conflict, we set forth a coping strategy that loosens the incompressibility constraint of damaged materials using the phase field degradation function. Such a scheme does not affect the incompressibility of intact materials, complying with the physical insights. With this idea in mind, the classic perturbed Lagrangian approach was modified, yielding a novel multi-field variational formulation for fracture problems. Then we derived the governing equations and outlined the concrete numerical implementation of the mixed formulation. The previously developed adaptive mesh deletion technology is also utilized to alleviate the mesh distortion problem in the 3D tearing test [49]. In addition, four representative numerical examples were chosen to demonstrate the validity and robustness of the proposed mixed framework. By comparing with the higher-order P2/P1 formulation, the superior performance of the classical Q1/P0 scheme in the modeling of nearly incompressible finite strain fractures was ascertained.

The rest of this paper is organized as follows. Section 2 emphatically uncovers the root cause of the incompressibility thwarting diffuse crack opening and conceives a countermeasure, thus laying the groundwork for the subsequently mixed formulation. In section 3, a variational form of the multi-field fracture problem is established after revisiting the perturbed Lagrangian approach. The governing equations in the strong form are also derived. In the aspect of numerical realization, Section 4 elaborates the discretization and linearization formulations, as well as some crucial numerical techniques. Section 5 demonstrates the validity of the proposed mixed framework via several representative numerical examples. Finally, Section 6 concludes this paper with the main contributions.

## 2. Phase field description for finite strain fracture with incompressible

### 2.1. Phase-field representation of cracks at finite strain

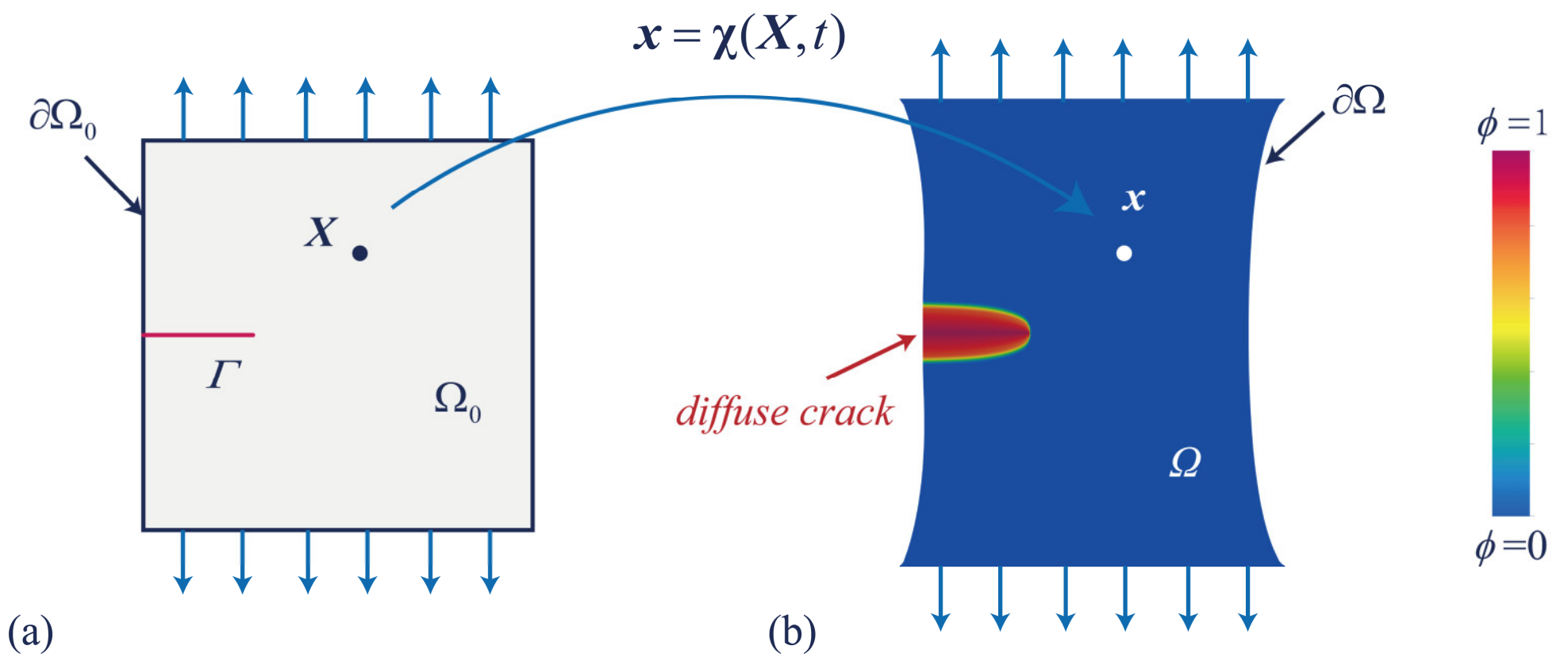

Fig. 1. (a) Illustration of a soft rectangular sheet containing a sharp crack $\Gamma$ in its undeformed configuration $\Omega_0$ . (b) Diffuse crack pattern depicted by phase field $\phi$ in the deformed configuration $\Omega$.

The highlight of phase field representation of cracks consists in smearing the sharp crack topology with a time-dependent scalar damage field $\phi(\boldsymbol{X},t)\in[0,1]$. To illustrate the evolution of diffuse cracks, we conceive a scenario where a soft rectangular sheet containing crack $\Gamma$ with initial configuration $\Omega_0 \subset \mathbb{R}^n$ and the boundary $\partial\Omega_0 \subset \mathbb{R}^{n-1}$ is stretched, as illustrated in Fig. 1(a). The separated crack surfaces that ought to open up parabolically in the new configuration $\Omega$ are now replaced by the phase field $\phi$ with a smooth transition from 0 (pristine state) to 1 (fully broken state), see Fig. 1(b). Let $\chi$ denote the deformation mapping from point $\boldsymbol{X}\in\Omega_0$ to $\boldsymbol{x}\in\Omega$ , then the displacement field is defined as

$$\boldsymbol{u} = \chi(\boldsymbol{X}) - \boldsymbol{X} \tag{1}$$

With this definition, the deformation gradient tensor $\boldsymbol{F}$, the Jacobian determinant $J$ and the right Cauchy–Green deformation tensor $\boldsymbol{C}$ are given by

$$\begin{aligned} \boldsymbol{F} &= \nabla_X \chi(\boldsymbol{X},t) = \boldsymbol{I} + \nabla_X \boldsymbol{u} \\ J &= \det(\boldsymbol{F}) \\ \boldsymbol{C} &= \boldsymbol{F}^T\boldsymbol{F} \end{aligned} \quad , \tag{2}$$

where $\boldsymbol{I}$ is the second-order identity tensor.

In the variational approach to fracture, the crack surface energy came by integrating the critical energy release of materials $G_C$ on the crack surfaces that can be approximated as

$$\Gamma_l(\phi) = \int_{\Gamma} G_c \mathrm{d}A \approx \int_{\Omega} G_c \gamma(\phi, \nabla\phi) \mathrm{d}V \tag{3}$$

where $\gamma$ is the crack surface density function devised complying with certain rules. Throughout most of the stories relating to PFF, a predominant format of $\gamma$ reads [20]

$$\gamma(\phi, \nabla\phi) = \frac{1}{4c_w}\left(\frac{w(\phi)}{l_0} + l_0 |\nabla\phi|^2\right) \text{ with } c_w = \int_0^1 \sqrt{w(\alpha)} d\alpha \tag{4}$$

Here, $l_0$ signifies the characteristic length absenting in the classical fracture theory, which dictates the width of the damage band of the smeared crack. In the present work, we intentionally opt for the following format

$$\gamma(\phi, \nabla\phi) = \frac{3}{8}\left(\frac{\phi}{l_0} + l_0 |\nabla\phi|^2\right), \left(w(\phi) = \phi, c_w = \frac{2}{3}\right) \tag{5}$$

The motivation for such selection emanates from the nature of the model that characterizes the pure elasticity state prior to the onset of damage.

Concerning the total potential energy functional at finite strain, a well-received potential energy functional take the form as [30]

$$\Pi_{u-\phi} = \int_{\Omega_0} g(\phi)\psi(\boldsymbol{C}) dV + \int_{\Omega_0} G_C \gamma(\phi, \nabla\phi) dV + \frac{\eta}{2}\int_{\Omega_0} \left(\frac{\partial\phi}{\partial t}\right)^2 dV - \Pi_{ext} \tag{6}$$

Herein, the degradation function takes $g(\phi) = (1-\phi)^2 + \varepsilon, (0 \le \varepsilon \ll 1)$ and $\psi(\boldsymbol{C})$ is the strain energy function. $\eta \ge 0$ (unit: $\mathrm{Nm}^{-2} \cdot \mathrm{s}$) is the artificial viscosity coefficient that regulates the energy dissipation rate, added for numerical robustness.

## 2.2. Incompressible hyperelastic model

Numerous hyperelastic materials can be regarded as nearly or truly incompressible, and the related constitutive models have also evolved into a big family. Herein, we postulate that the strain energy function $\psi(\boldsymbol{C})$ follows a general decomposition form [50]

$$\psi(\boldsymbol{C}) = \tilde{\psi}(\boldsymbol{C}) + \underbrace{\kappa U(J)}_{\psi^{vol}} \tag{7}$$

where $\kappa$ is the bulk modulus, and the volumetric part $\psi^{vol} = \kappa U(J)$ vanishes in the incompressible limit. Concerning $\tilde{\psi}(\boldsymbol{C})$, let us focus on a sophisticated model with a phenomenological foundation developed by Ogden [51], which takes the form

$$\tilde{\psi}(\lambda_1, \lambda_2, \lambda_3) = \sum_{a=1}^{N} \frac{\mu_a}{m_a}\left(\lambda_1^{m_a} + \lambda_2^{m_a} + \lambda_3^{m_a} - 3\right) \tag{8}$$

Note that $\lambda_i$ $(i=1,2,3)$ stand for the principal stretches, ruled by the spectral decomposition of the right Cauchy–Green deformation tensor $\boldsymbol{C}$

$$\boldsymbol{C}=\sum_{a=1}^{3}\lambda_a^2\boldsymbol{N_a}\otimes\boldsymbol{N_a} \tag{9}$$

In Eq. 8, $\mu_a$ and $m_a$ are Ogden parameters, and $N$ takes a positive integer. By setting $N=2,\ m_1=2,\ m_2=-2$ and keeping the incompressible constraint $J=\lambda_1\lambda_2\lambda_3=1$ in mind, we find that Eq.8 degenerates to the widely used Mooney-Rivlin model

$$\begin{aligned}\tilde{\psi}_{MR}&=c_1\left(\lambda_1^2+\lambda_2^2+\lambda_3^2-3\right)+c_2\left(\lambda_1^{-2}+\lambda_2^{-2}+\lambda_3^{-2}-3\right)\\&=c_1\left(I_1-3\right)+c_2\left(I_2-3\right)\end{aligned} \tag{10}$$

Here, $I_1$ and $I_2$ are the first and second principal invariants, respectively. The conventional shear modulus $\mu$ has the value of $2(c_1+c_2)$ . Furthermore, a formulaically simpler neo-Hookean model

$$\tilde{\psi}_{NH}=c_1\left(\lambda_1^2+\lambda_2^2+\lambda_3^2-3\right)=c_1\left(I_1-3\right) \tag{11}$$

with $c_1=\mu/2$ can be attained by setting $N=1,\ m_1=2$ .

Considering the mathematical form of the strain-energy function $\tilde{\psi}$ expressed in the principal stretches $\lambda_1,\lambda_2,\lambda_3$, the second Piola–Kirchhoff stress (PK2) is derived as

$$\tilde{\boldsymbol{S}}=\sum_{a=1}^{3}\underbrace{\frac{1}{\lambda_a}\frac{\partial\tilde{\psi}}{\partial\lambda_a}}_{\tilde{S}_a}\boldsymbol{N_a}\otimes\boldsymbol{N_a}, \tag{12}$$

the first Piola–Kirchhoff stress (PK1) can thus be written straightforwardly via $\tilde{\boldsymbol{P}}=\boldsymbol{F}\cdot\tilde{\boldsymbol{S}}$. To facilitate the subsequent numerical treatment, we also derive the fourth-order elastic tensor $\tilde{\mathbb{C}}$ in spectral form

$$\tilde{\mathbb{C}}=\sum_{a,b=1}^{3}\frac{1}{\lambda_b}\frac{\partial\tilde{S}_a}{\partial\lambda_b}\mathbb{Q}_{ab}+\sum_{\substack{a,b=1,\ a\neq b\\ \lambda_a\neq\lambda_b}}^{3}\frac{\tilde{S}_b-\tilde{S}_a}{\lambda_b^2-\lambda_a^2}\mathbb{G}_{ab}+\sum_{\substack{a,b=1,\ a\neq b\\ \lambda_a=\lambda_b}}^{3}\frac{\partial\tilde{S}_b}{\partial\lambda_b^2}-\frac{\partial\tilde{S}_a}{\partial\lambda_b^2}\mathbb{G}_{ab} \tag{13}$$

with

$$\begin{aligned}\mathbb{Q}_{ab}&=\boldsymbol{N_a}\otimes\boldsymbol{N_a}\otimes\boldsymbol{N_b}\otimes\boldsymbol{N_b}\\ \mathbb{G}_{ab}&=\boldsymbol{N_a}\otimes\boldsymbol{N_b}\otimes\boldsymbol{N_a}\otimes\boldsymbol{N_b}+\boldsymbol{N_a}\otimes\boldsymbol{N_b}\otimes\boldsymbol{N_b}\otimes\boldsymbol{N_a}\end{aligned} \tag{14}$$

**Remark.** An alternative decomposition format of Eq. 7 takes the form as [52]

$$\psi\left(\boldsymbol{C}\right)=\psi^{iso}\left(\bar{\boldsymbol{C}}\right)+\psi^{vol}\left(J\right), \tag{15}$$

where $\psi^{iso}$ and $\psi^{vol}$ are the decoupled isochoric and volumetric parts, respectively. Accordingly, the PK2 stress and elastic tensors can be written as

$$\boldsymbol{S}=\boldsymbol{S}^{iso}+\boldsymbol{S}^{vol}=J^{-2/3}\underbrace{\left(\boldsymbol{I}-\frac{1}{3}\boldsymbol{C}^{-1}\otimes\boldsymbol{C}^{-1}\right)}_{\mathbb{P}}:2\frac{\partial\psi^{iso}\left(\bar{\boldsymbol{C}}\right)}{\partial\bar{\boldsymbol{C}}}+Jp\boldsymbol{C}^{-1} \tag{16}$$

and

$$
\begin{aligned}
\mathbb{C} &= \mathbb{C}^{iso} + \mathbb{C}^{vol} \\
&= \mathbb{P} : 4J^{-4/3} \frac{\partial^2 \psi^{iso}\left(\bar{\boldsymbol{C}}\right)}{\partial \bar{\boldsymbol{C}} \partial \bar{\boldsymbol{C}}} : \mathbb{P} \\
&+ \frac{2}{3}\left( J^{-2/3} \frac{\partial \psi^{iso}\left(\bar{\boldsymbol{C}}\right)}{\partial \bar{\boldsymbol{C}}} : \boldsymbol{C} \right)\left( \boldsymbol{C}^{-1} \odot \boldsymbol{C}^{-1} - \frac{1}{3} \boldsymbol{C}^{-1} \otimes \boldsymbol{C}^{-1} \right) \\
&- \frac{2}{3}\left( \boldsymbol{C}^{-1} \otimes \boldsymbol{S}^{iso} - \boldsymbol{S}^{iso} \otimes \boldsymbol{C}^{-1} \right)
\end{aligned}
\tag{17}
$$

Note that, the above formula is appealing in the research of nearly incompressible hyperelastic materials, whereas our pre-tests indicate that it lacks numerical robustness in the phase field modeling of fracture coupled with incompressibility. Therefore, this decomposition form is deprecated in this work.

## 2.3. Incompressibility - a barrier for diffuse crack opening

The diffusive depiction of sharp cracks circumvents the laborious tracking of discontinuous surfaces, which is the principal cause that PFM is favored in fracture modeling. However, this advantage may turn into a barrier in some special circumstances, such as finite strain fracture of incompressible hyperelastic materials. Revisiting Fig. 1, we remark that the physically material-free crack opening zone is still integrated into the computational domain, which inevitably generates volume expansion compared to discrete crack configurations. In this context, enforcing the incompressible constraint will frustrate the parabolic opening of the crack, and even result in the collapse of the numerical solution (see **Appendix A** for an intuitive presentation). Given this intrinsic conflict, simply seeking an alternative solution algorithm is doomed to be futile. Enlightened by the idea of attenuating strain energy via multiplying the degradation function $g(\phi)$ in the PFF modeling, we loosen the incompressible constraint in the same manner. Now, using $\vartheta_0 = \kappa_0 / \mu_0$ to measure the incompressible level, we assume that

$$\vartheta(\phi) = f(\phi)\vartheta_0 = f(\phi)\frac{\kappa_0}{\mu_0} \tag{18}$$

where $f\left(\phi\right)$ satisfies

$$f\left(0\right) = 0, f\left(1\right) = 1 \tag{19}$$

This definition indicates that formula 18 only relieves the incompressible constraint of the damaged phase, and yet does not loosen the binding for non-damaged materials. Considering that the damaged material within the crack opening outline has no substantial contribution to the physical entity, such a treatment is sound in the limit of $\Gamma$-convergence [53]. Moreover, due to the moduli $\kappa$ and $\mu$ have been degenerated in the pristine phase field formulation, Eq.18 can be rephrased as

$$\vartheta(\phi) = f(\phi)\frac{g(\phi)\kappa_0}{g(\phi)\mu_0} = \frac{f(\phi)g(\phi)\kappa_0}{g(\phi)\mu_0} \tag{20}$$

For simplicity, we let $f(\phi) = g(\phi)$, thus yielding

$$\vartheta(\phi) = \frac{g(\phi)^2 \kappa_0}{g(\phi)\mu} \tag{21}$$

This expression implies that the bulk modulus $\kappa$ requires to be damaged faster than the shear modulus $\mu$. For the truly incompressible hyperelastic materials, however, the above countermeasure fails as the bulk modulus $\kappa$ approaches infinity. Consequently, slightly compressible (quasi-incompressible) is introduced in the modeling, arriving at the same destination as the subsequent perturbed Lagrangian approach.

## 3. Mixed formulation

### 3.1. Perturbed Lagrangian approach

To date, the classical Lagrangian multiplier approach is the sole workable technique for fully incompressible cases [54]. Without involving phase field damage, the potential energy functional of the Lagrangian multiplier approach is written as

$$\Pi_{LagMul}(\boldsymbol{u}, p) = \int_{\Omega_0} \left[\tilde{\psi}(\boldsymbol{C}) + p(J-1)\right] dV - \Pi_{ext} \tag{22}$$

with the external potential energy ($\Pi_{ext}$) comprised of body and surface loads

$$\Pi_{ext} = \int_{\Omega_0} \bar{\boldsymbol{b}}_0 \cdot \boldsymbol{u} dV + \int_{\partial\Omega_0} \bar{\boldsymbol{t}}_n \cdot \boldsymbol{u} dA \tag{23}$$

Herein, the scalar $p$ acts as the Lagrangian multiplier, which can be recognized as an independent variable representing hydrostatic pressure. Nonetheless, this formulation yields a saddle-point problem that is prone to numerical difficulties. To address the numerical issue, a perturbed Lagrangian approach is proposed [45], and the corresponding energy functional is reformatted as

$$\Pi_{PerMul}(\boldsymbol{u}, p) = \int_{\Omega_0} \left( \tilde{\psi}(\boldsymbol{C}) + p\left( J - 1 - \frac{p}{2\kappa} \right) \right) dV - \Pi_{ext} \tag{24}$$

Note that, such an approach is compatible with nearly and fully incompressible materials, thanks to slightly compressible modifications. In the incompressible limit, i.e., $\kappa \to +\infty$, the Lagrangian multiplier approach (Eq. 22) can be recovered from Eq. 24. It is worth mentioning that Kadapa and Hossain have proposed a unified framework that further generalized the perturbed Lagrangian method to incorporate compressible materials [55]. However, we do not intend to embrace this formulation, allowing for the keynote of the current work.

### 3.2. Variational form for multi-field fracture problem

In the previous displacement-phase field ($\boldsymbol{u}-\phi$) solution scheme, the mechanics sub-problem only solves a single displacement field. We remark that such an approach is incapable of modeling nearly and truly incompressible materials. Thus the functional $\Pi_{u-\phi}$ requires to be amended in the perturbed Lagrangian approach. A straightforward extension is to multiply the perturbed Lagrangian term by $g(\phi)$, resulting in

$$\begin{aligned}\Pi^{test}_{\boldsymbol{u}-p-\phi} &= \int_{\Omega_0} g(\phi)\left(\tilde{\psi}(\boldsymbol{C}) + p\left(J-1-\frac{p}{2\kappa}\right)\right)dV + \int_{\Omega_0} G_C\gamma(\phi,\nabla\phi)dV \\ &+\frac{\eta}{2}\int_{\Omega_0}\left(\frac{\partial\phi}{\partial t}\right)^2 dV - \Pi_{ext}\end{aligned} \tag{25}$$

This format excels in modeling incompressible (quasi-incompressible) brittle fracture [56]. however, it fails in the finite deformation context due to the forced incompressibility of the damaged materials (see **Appendix A**). For this reason, guided by the ideas illuminated in sub-section 2.3, we liberate the incompressibility constraints of the damaged phase, namely $\kappa = \left(g(\phi)\right)^2 \kappa_0$, to eliminate the barriers to crack opening. Assuming the damaged volumetric energy function

$$\psi_{vol} = \frac{1}{2}\left(g(\phi)\right)^2 \kappa_0 \left(J-1\right)^2 \tag{26}$$

thus, the pressure takes the form

$$p = \frac{\partial \psi_{vol}}{\partial J} = \left(g(\phi)\right)^2 \kappa_0 \left(J-1\right) \tag{27}$$

Remarkably, this expression can be identified as a constraint, i.e., $\left(g(\phi)\right)^2\left(J-1\right) - \frac{p}{\kappa_0} = 0$, enforced via a Lagrange multiplier $p$. In this manner, we now arrive at mixed energy functional for the multi-field fracture problem

$$\begin{aligned}\Pi_{u-p-\phi} &= \int_{\Omega_0} g(\phi)\tilde{\psi}(\boldsymbol{C})dV + \int_{\Omega_0} p\left(\left(g(\phi)\right)^2 (J-1) - \frac{p}{2\kappa_0}\right)dV + \\ &\int_{\Omega_0} G_C\gamma(\phi,\nabla\phi)dV + \frac{\eta}{2}\int_{\Omega_0}\left(\frac{\partial\phi}{\partial t}\right)^2 dV - \Pi_{ext}\end{aligned} \tag{28}$$

Taking the variational of Eq. 28 in terms of $\boldsymbol{u}, p, \phi$, we obtain

$$\delta\Pi_{u-p-\phi}=\int_{\Omega_0}\left(2g(\phi)\frac{\partial\tilde{\psi}(\boldsymbol{C})}{\partial\boldsymbol{C}}\boldsymbol{F}^T:\nabla\delta\boldsymbol{u}+\frac{\partial g(\phi)}{\partial\phi}\tilde{\psi}(\boldsymbol{C})\delta\phi\right)dV$$
$$+\int_{\Omega_0}\left(g(\phi)\right)^2 pJ\boldsymbol{C}^{-1}\boldsymbol{F}^T:\nabla\delta\boldsymbol{u}dV+\int_{\Omega_0}2g(\phi)\frac{\partial g(\phi)}{\partial\phi}p(J-1)\delta\phi dV$$
$$+\int_{\Omega_0}\left(\left(g(\phi)\right)^2(J-1)-\frac{p}{\kappa_0}\right)\delta p dV+\int_{\Omega_0}\frac{3}{8}G_C\left(\frac{1}{l_0}\delta\phi+2l_0\nabla\phi\cdot\nabla\delta\phi\right)dV \quad (29)$$
$$+\int_{\partial\Omega_0}\frac{3}{4}G_c l_0\nabla_X\phi\cdot\boldsymbol{n}_0\delta\phi\mathrm{d}A+\eta\int_{\Omega_0}\frac{\partial\phi}{\partial t}\delta\phi dV-\Pi_{ext}(\delta\boldsymbol{u})$$

After taking the integration by parts, Eq.29 is reformulated as

$$\delta\Pi_{u-p-\phi}=-\int_{\Omega_0}\left(\nabla_X\cdot\tilde{\boldsymbol{P}}\right)\cdot\delta\boldsymbol{u}\mathrm{d}V-\int_{\Omega_0}\left(\nabla_X\cdot\boldsymbol{P}^{vol}\right)\cdot\delta\boldsymbol{u}dV$$
$$+\int_{\partial\Omega_0^N}\left((\tilde{\boldsymbol{P}}+\boldsymbol{P}^{vol})\cdot\boldsymbol{n}\right)\cdot\delta\boldsymbol{u}\mathrm{d}A+\int_{\Omega_0}\left(\left(g(\phi)\right)^2(J-1)-\frac{p}{\kappa_0}\right)\delta p dV$$
$$+\int_{\Omega_0}\left(\frac{\partial g(\phi)}{\partial\phi}\tilde{\psi}(\boldsymbol{C})+\frac{3}{8}\left(\frac{G_c}{l_0}-2G_c l_0\nabla^2\phi\right)+2g(\phi)\frac{\partial g(\phi)}{\partial\phi}p(J-1)+\eta\dot{\phi}\right)\delta\phi\mathrm{d}V \quad (30)$$
$$+\int_{\partial\Omega_0}\frac{3}{4}G_c l_0\nabla_X\phi\cdot\boldsymbol{n}_0\delta\phi\mathrm{d}A-\int_{\Omega_0}\bar{\boldsymbol{b}}_0\cdot\delta\boldsymbol{u}dV+\int_{\partial\Omega_0}\bar{\boldsymbol{t}}_0\cdot\delta\boldsymbol{u}dA$$

with the nominal stress

$$\tilde{\boldsymbol{P}}=2g(\phi)\frac{\partial\tilde{\psi}(\boldsymbol{C})}{\partial\boldsymbol{C}}\boldsymbol{F}^T \quad (31)$$

and

$$\boldsymbol{P}^{vol}=\left(g(\phi)\right)^2 Jp\boldsymbol{C}^{-1}\boldsymbol{F}^T \quad (32)$$

### 3.3. Governing equations

On account of the independence of the variational $\delta\boldsymbol{u}$, $\delta p$, $\delta\phi$, the equilibrium condition of formula 30 is to satisfy the following three equations

$$-\int_{\Omega_0}\left(\nabla_X\cdot\left(\tilde{\boldsymbol{P}}+\boldsymbol{P}^{vol}\right)-\bar{\boldsymbol{b}}_0\right)\cdot\delta\boldsymbol{u}\mathrm{d}V+\int_{\partial\Omega_0^N}\left((\tilde{\boldsymbol{P}}+\boldsymbol{P}^{vol})\cdot\boldsymbol{n}-\bar{\boldsymbol{t}}_0\right)\cdot\delta\boldsymbol{u}\mathrm{d}A=0 \quad (33)$$

$$\int_{\Omega_0}\left(\left(g(\phi)\right)^2(J-1)-\frac{p}{\kappa_0}\right)\delta p dV=0 \quad (34)$$

$$\int_{\Omega_0}\left(\frac{\partial g(\phi)}{\partial\phi}\tilde{\psi}(\boldsymbol{C})+\frac{3}{8}\left(\frac{G_c}{l_0}-2G_c l_0\nabla^2\phi\right)+2g(\phi)\frac{\partial g(\phi)}{\partial\phi}p(J-1)+\eta\dot{\phi}\right)\delta\phi\mathrm{d}V$$
$$+\int_{\partial\Omega_0}\frac{3}{4}G_c l_0\nabla_X\phi\cdot\boldsymbol{n}_0\delta\phi\mathrm{d}A=0 \quad (35)$$

Notwithstanding the strong-form governing equations are unessential for numerical treatment, we still derived them for the sake of mathematical succinct, as follows

$$
\begin{cases}
\nabla_X \cdot (\tilde{\boldsymbol{P}} + \boldsymbol{P}^{vol}) + \bar{\boldsymbol{b}}_0 = 0 \\
\left(g(\phi)\right)^2 (J-1) - \dfrac{p}{\kappa_0} = 0 \\
\dfrac{\partial g(\phi)}{\partial \phi}\tilde{\psi}(\boldsymbol{C}) + \dfrac{3}{8}(\dfrac{G_c}{l_0} - 2G_c l_0 \nabla^2 \phi) + 2g(\phi)\dfrac{\partial g(\phi)}{\partial \phi} p(J-1) + \eta\dot{\phi} = 0 \\
(\tilde{\boldsymbol{P}} + \boldsymbol{P}^{vol}) \cdot \boldsymbol{n} = \bar{\boldsymbol{t}}_0 \qquad \text{at } \partial\Omega_0^N \\
\nabla_X \phi \cdot \boldsymbol{n}_0 = 0 \qquad \text{at } \partial\Omega_0
\end{cases}
\tag{36}
$$

Note that, Eq. 27, which serves to relax the incompressible constraint of the damage phase, is strictly derived, see the second equality of Eq. 36.

## 4. Numerical aspects

### 4.1. Weak form of mixed formulation

The standardized finite element method (FEM) solving derives the weak form of the governing equation as per the Galerkin weighted residual method at first. Let $\boldsymbol{w}^{\boldsymbol{u}}$, $w^p$, $w^\phi$ denote the test functions for the mechanical response, pressure, and phase field equations, respectively. Multiply the strong form (Eq. 36) by these test functions and perform integration by parts, the weak form of mixed formulation in the residual form can then be derived as

$$
\boldsymbol{R}_{\boldsymbol{u}} = \int_{\Omega_0} \left(\tilde{\boldsymbol{P}} + \boldsymbol{P}^{vol}\right) : \nabla_X \boldsymbol{w}^{\boldsymbol{u}} \mathrm{d}V - \int_{\Omega_0} \bar{\boldsymbol{b}}_0 \cdot \boldsymbol{w}^{\boldsymbol{u}} \mathrm{d}V - \int_{\Gamma_N} \bar{\boldsymbol{t}}_0 \cdot \boldsymbol{w}^{\boldsymbol{u}} \mathrm{d}A = 0 \text{ }, \forall \boldsymbol{w}^{\boldsymbol{u}} \in \mathcal{W}_{\boldsymbol{u}} \tag{37}
$$

$$
\boldsymbol{R}_p = \int_{\Omega_0} \left( \left( (1-\phi)^4 + \varepsilon \right)(J-1) - \frac{p}{\kappa_0} \right) w^p \mathrm{d}V = 0 \text{ }, \forall w^p \in \mathcal{W}_p \tag{38}
$$

$$
\begin{aligned}
\boldsymbol{R}_\phi = & \int_{\Omega_0} \left( -2(1\text{-}\phi)\tilde{\psi}(\boldsymbol{C}) - 4(1\text{-}\phi)^3 p(J-1) \right) w^\phi \mathrm{d}V \\
& + \int_{\Omega_0} \left( \frac{3}{8}\frac{G_c}{l_0}\left(2l_0^2 \nabla_X \phi \cdot \nabla_X w^\phi + w^\phi\right) + \eta\dot{\phi} w^\phi \right) \mathrm{d}V = 0 \qquad \forall w^\phi \in \mathcal{W}_\phi
\end{aligned}
\tag{39}
$$

Herein, $\mathcal{W}_{\boldsymbol{u}}$, $\mathcal{W}_p$ and $\mathcal{W}_\phi$ are trial function spaces for displacement field $\boldsymbol{u}$, pressure field $p$ and phase-field $\phi$.

### 4.2. FEM discretization

To assure the stability of mixed finite element formulation, the collocation of discretization elements should be selected elaborately. Now, we exploit the well-known Q1/P0 element [38], as well as the P2/P1 element (a member of the Taylor-Hood element family [57]) with *inf-sup* stability for discretization, and the corresponding configuration and node layout of these two elements are illustrated in Fig. 2. In the Q1/P0 formulation, displacement $\boldsymbol{u}$ and phase-field $\phi$ are discretized by linear quadrilateral element (Q4) for the 2D case or linear hexahedral elements (H8) for the 3D case, while the pressure $p$ retains constant within an element. Concerning the P2/P1 formulation, $\boldsymbol{u}$ is discretized using quadratic 6-noded triangular (10-noded tetrahedral) elements, whilst $p$ and $\phi$ are discretized with linear 3-noded (4-noded) ones.

As a demonstration, the Q1/P0 element is applied for the discretization in space, thus the fundamental field variables $\boldsymbol{u}$, $p$, $\phi$ and their functions $\boldsymbol{w}^{\boldsymbol{u}}$, $w^p$, $w^\phi$ are approximated as

$$\begin{aligned} &\boldsymbol{u}=\sum_{i=1}^{m}\boldsymbol{N}_i^{\boldsymbol{u}}\boldsymbol{u}_i,\ \ p=\sum_{i=1}^{n}N_i^p p_i,\ \ \phi=\sum_{i=1}^{m}N_i^\phi\phi_i \\ &\boldsymbol{w}^{\boldsymbol{u}}=\sum_{i=1}^{m}\boldsymbol{N}_i^{\boldsymbol{u}}\boldsymbol{w}_i^{\boldsymbol{u}},\ w^p=\sum_{i=1}^{n}N_i^p w_i^p,\ w^\phi=\sum_{i=1}^{m}N_i^\phi w_i^\phi \end{aligned} \tag{40}$$

with

$$\boldsymbol{N}_i^{\boldsymbol{u}}=\begin{bmatrix} N_i^u & 0 \\ 0 & N_i^u \end{bmatrix} \tag{41}$$

Herein, $N_i^u=N_i^\phi$ are the standard basis functions for linear quadrilateral (hexahedral) elements, while $N_i^p=1$. As a result, the nodal number in each element $m=4$ (8) and $n=1$.

In regard of the time-related dissipation term $\eta\dot{\phi}$, a concise backward Euler difference method is adopted for the temporal discretization

$$\dot{\phi}=\frac{\phi_{n+1}-\phi_n}{\Delta t} \tag{42}$$

where the subscript n+1 tags the current time $t_{n+1}$. The time step $\Delta t$ used in this work is set to the order of $10^{-5}\mathrm{s}\sim10^{-7}\mathrm{s}$.

**Remark.** Albeit suffering some complaints due to the lack of *inf-sup* stability, the classical Q1/P0 formulation works very robustly in most computing scenarios involving (nearly) incompressibility, containing the current issues. As for the higher-order elements, we exclusively chose the P2/P1 element because of the enormously fine mesh demand in PFF modeling. Other alternative higher-order elements such as BQ2/BQ1 are also available [58], and yet the computational cost is daunting, especially in 3-D phase field modeling of fracture.

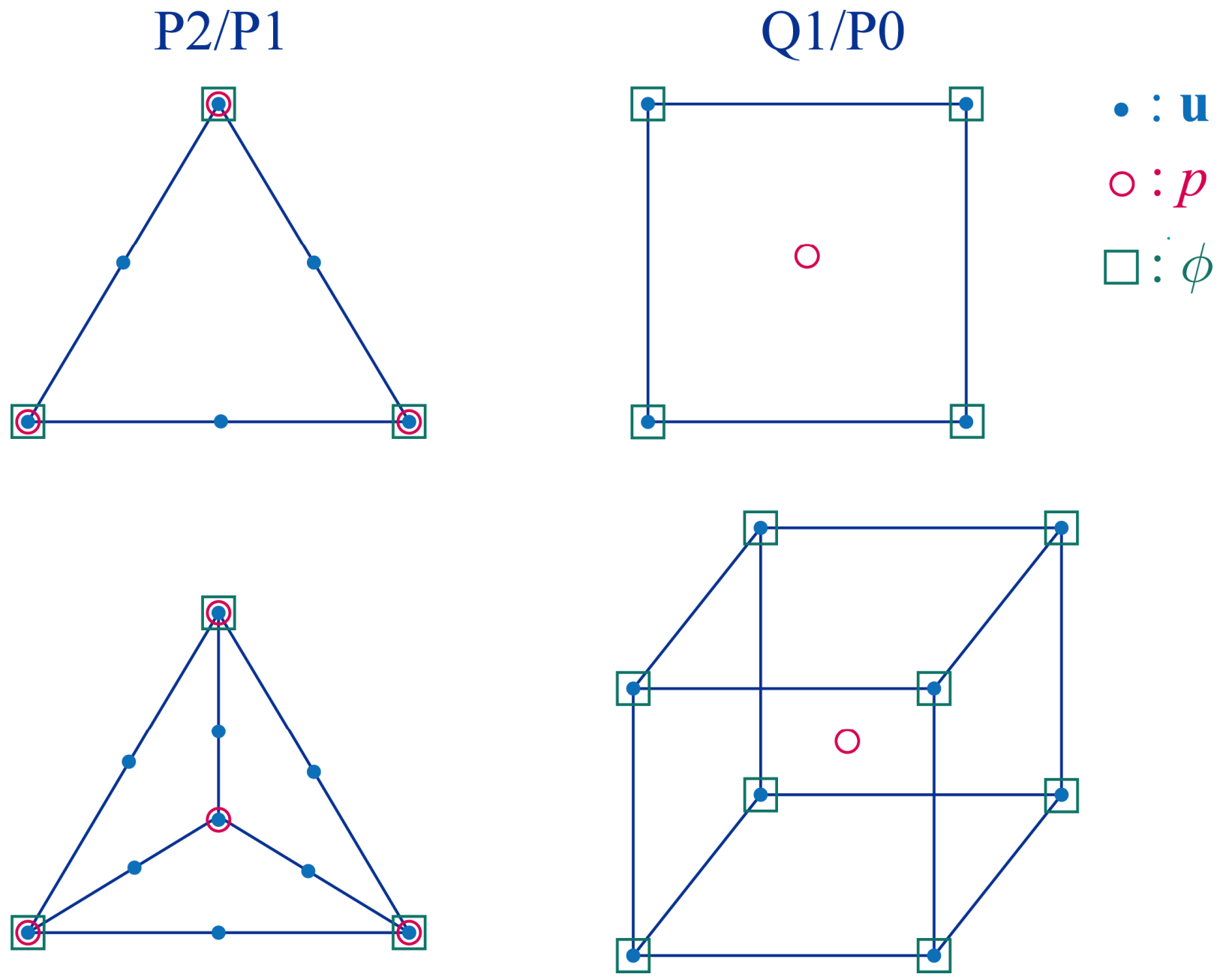

Fig. 2. Shape and node arrangement of Q1/P0 and P2/P1 elements in 2D and 3D cases. Symbol •, ○ and □ indicate displacement $\boldsymbol{u}$, pressure $p$ and phase field $\phi$ DOFs, respectively.

## 4.3. Linearization

To solve the coupled Eqs 37-38, we linearize them using the Newton-Raphson method, as follows

$$\begin{bmatrix} \boldsymbol{K}^e_{\boldsymbol{uu}} & \boldsymbol{K}^e_{\boldsymbol{u}p} & \boldsymbol{K}^e_{\boldsymbol{u}\phi} \\ \boldsymbol{K}^e_{p\boldsymbol{u}} & \boldsymbol{K}^e_{pp} & \boldsymbol{K}^e_{p\phi} \\ \boldsymbol{K}^e_{\phi\boldsymbol{u}} & \boldsymbol{K}^e_{\phi p} & \boldsymbol{K}^e_{\phi\phi} \end{bmatrix} \begin{Bmatrix} \Delta\boldsymbol{u} \\ \Delta p \\ \Delta\phi \end{Bmatrix} = \begin{Bmatrix} \boldsymbol{f}^e_{\boldsymbol{u}} \\ \boldsymbol{f}^e_p \\ \boldsymbol{f}^e_\phi \end{Bmatrix} \tag{43}$$

Here, Eq. 43 involves the three-field coupling of $(\boldsymbol{u}, p, \phi)$, and its monolithic solution is generally challenging. Imitating the tried-and-tested staggered algorithm, we decoupled the displacement-pressure field $(\boldsymbol{u}, p)$ and the phase field $\phi$ such that

$$\boldsymbol{K}^e_{\phi\boldsymbol{u}} = \boldsymbol{K}^e_{\phi p} = \boldsymbol{K}^e_{\boldsymbol{u}\phi} = \boldsymbol{K}^e_{p\phi} = \boldsymbol{0} \tag{44}$$

Thus, Eq. 43 is recast as

$$\begin{cases} \boldsymbol{K}^e_{\boldsymbol{uu}} \cdot \Delta\boldsymbol{u} + \boldsymbol{K}^e_{\boldsymbol{u}p} \cdot \Delta p = \boldsymbol{f}^e_{\boldsymbol{u}} \\ \boldsymbol{K}^e_{p\boldsymbol{u}} \cdot \Delta\boldsymbol{u} + \boldsymbol{K}^e_{pp} \cdot \Delta p = \boldsymbol{f}^e_p \\ \boldsymbol{K}^e_{\phi\phi} \cdot \Delta\phi = \boldsymbol{f}^e_\phi \end{cases} \tag{45}$$

By introducing the discretization defined by Eq. 40, the residual vectors $\boldsymbol{f}^e_{\boldsymbol{u}}$, $\boldsymbol{f}^e_p$ and $\boldsymbol{f}^e_\phi$ at the element level take on the form

$$\begin{cases} \boldsymbol{f}_{\boldsymbol{u}}^{e} = -\int_{\Omega_0} \boldsymbol{B}_X^T \left\{ \tilde{\boldsymbol{S}} + \boldsymbol{S}^{vol} \right\} \mathrm{d}V + \int_{\Omega_0} N^{\boldsymbol{u}} \bar{\boldsymbol{b}} \mathrm{d}V + \int_{\Gamma_N} N^{\boldsymbol{u}} \bar{\boldsymbol{t}} \mathrm{d}A \\ \boldsymbol{f}_p^e = \int_{\Omega_0} \left( \left( (1-\phi)^4 + \varepsilon \right)(J-1) - \dfrac{p}{\kappa_0} \right) \mathrm{d}V \\ \boldsymbol{f}_\phi^s = \int_{\Omega_0} \left( 2(1-\phi_{n+1}) N^\phi \tilde{\psi}(\boldsymbol{C}) + 4(1\text{-}\phi)^3 N^\phi p (J-1) \right) \mathrm{d}V \\ -\int_{\Omega_0} \left( \dfrac{3}{4} \boldsymbol{B}_\phi^T G_c l_0 \boldsymbol{B}_\phi \phi_{n+1} + \dfrac{3}{8} G_c l_0 N^\phi + \dfrac{\eta}{\Delta t} N^\phi \left( \phi_{n+1} - \phi_n \right) \right) dV \end{cases} \tag{46}$$

Going by the matrix derivation rule, the tangent stiffness matrices thereof $\boldsymbol{K}_{\boldsymbol{uu}}^e, \boldsymbol{K}_{\boldsymbol{up}}^e$, $\boldsymbol{K}_{\boldsymbol{pu}}^e$, $\boldsymbol{K}_{pp}^e$ and $\boldsymbol{K}_{\phi\phi}^e$ are derived as

$$\begin{cases} \boldsymbol{K}_{\boldsymbol{uu}}^e = \dfrac{\partial R_{\boldsymbol{u}}^e}{\partial \boldsymbol{u}} = \int_{\Omega_0} \boldsymbol{B}_X^T \left( \hat{\mathbb{C}} + \hat{\mathbb{C}}^{vol} \right) \boldsymbol{B}_X \mathrm{d}V + \int_{\Omega_k} \mathcal{B}^T \left( \hat{\boldsymbol{S}} + \hat{\boldsymbol{S}}^{vol} \right) \mathcal{B} \mathrm{d}V \\ \boldsymbol{K}_{\boldsymbol{up}}^e = \boldsymbol{K}_{\boldsymbol{pu}}^{e\ T} = \dfrac{\partial R_{\boldsymbol{u}}^e}{\partial \phi} = \int_{\Omega_0} \boldsymbol{B}_X^T \left( (1-\phi)^4 + \varepsilon \right) J \boldsymbol{C}^{-1} \mathrm{d}V \\ \boldsymbol{K}_{pp}^e = \dfrac{\partial R_p^e}{\partial p} = -\int_{\Omega_0} \dfrac{1}{\kappa_0} \mathrm{d}V \\ \boldsymbol{K}_{\phi\phi}^e = \int_{\Omega_0} \left[ 2N^\phi \tilde{\psi}(\boldsymbol{C}) N^{\phi^T} + 12(1-\phi)^2 N^{\phi^T} p(J-1) + \dfrac{3}{4} \boldsymbol{B}_\phi^T G_c l_0 \boldsymbol{B}_\phi + \dfrac{\eta}{\Delta t} N^\phi N^{\phi^T} \right] \mathrm{d}V \end{cases} \tag{47}$$

Herein, the damaged PK2 stress $\hat{\boldsymbol{S}}$ and elastic tensor $\hat{\mathbb{C}}$ read

$$\hat{\boldsymbol{S}} = \left( (1-\phi)^2 + \varepsilon \right) \tilde{\boldsymbol{S}} \tag{48}$$

and

$$\hat{\mathbb{C}} = \left( (1-\phi)^2 + \varepsilon \right) \tilde{\mathbb{C}} \tag{49}$$

respectively, while the pure volumetric parts are given by

$$\boldsymbol{S}^{vol} = \boldsymbol{P}^{vol} \cdot \boldsymbol{F}^{-T} = \left( (1-\phi)^4 + \varepsilon \right) Jp\boldsymbol{C}^{-1} \tag{50}$$

$$\mathbb{C}^{vol} = 2\frac{\partial \boldsymbol{S}^{vol}}{\partial \boldsymbol{C}} = \left( (1-\phi)^4 + \varepsilon \right) \left( pJ\, \mathrm{C}_{ij}^{-1} \mathrm{C}_{kl}^{-1} - pJ \left( \mathrm{C}_{il}^{-1} \mathrm{C}_{kj}^{-1} + \mathrm{C}_{ik}^{-1} \mathrm{C}_{jl}^{-1} \right) \right) \tag{51}$$

Besides, the gradient matrices $\boldsymbol{B}_X$ and $\mathcal{B}$ are defined by

$$\boldsymbol{B}_X = \begin{bmatrix} N_{i,X} F_{11} & N_{i,X} F_{21} & N_{i,X} F_{31} \\ N_{i,X} F_{12} & N_{i,Y} F_{22} & N_{i,Y} F_{32} \\ N_{i,Z} F_{13} & N_{i,Z} F_{23} & N_{i,Z} F_{33} \\ N_{i,Y} F_{13} + N_{i,Z} F_{12} & N_{i,Y} F_{23} + N_{i,Z} F_{22} & N_{i,Y} F_{33} + N_{i,Y} F_{32} \\ N_{i,X} F_{13} + N_{i,Z} F_{11} & N_{i,X} F_{23} + N_{i,Z} F_{21} & N_{i,X} F_{33} + N_{i,Z} F_{31} \\ N_{i,X} F_{12} + N_{i,Y} F_{11} & N_{i,X} F_{22} + N_{i,Y} F_{21} & N_{i,X} F_{32} + N_{i,Y} F_{31} \end{bmatrix} \tag{52}$$

and

$$\mathcal{B}^T = \begin{bmatrix} N_{i,X} & N_{i,Y} & N_{i,Z} & 0 & 0 & 0 & 0 & 0 & 0 \\ 0 & 0 & 0 & N_{i,X} & N_{i,Y} & N_{i,Z} & 0 & 0 & 0 \\ 0 & 0 & 0 & 0 & 0 & 0 & N_{i,X} & N_{i,Y} & N_{i,Z} \end{bmatrix} \tag{53}$$

where $N_{i,X}$ symbolizes $\partial N_i / \partial X$ .

### 4.4. Enforced constraints

As far as the topic of this article is concerned, we touch upon two categories of constraints, i.e., the incompressible constraints and the irreversible constraints on crack growth. The basic tenet for disposing of the former is to degenerate the incompressibility of damaged materials, which has been elaborated earlier and will not be reiterated here. Referring to the irreversible constraints

$$\dot{\phi} = \frac{\phi_{n+1} - \phi_{n+1}}{\Delta t} \geq 0 \tag{54}$$

which is imposed by using an active set method that can be well embedded in Newton–Raphson iterations [27, 49]. For the sake of clarity, we narrate this method in the form of pseudocode, as stated in Algorithm 1.

**Algorithm 1** Active set method.

1. Define $\mathcal{A} = \varnothing, \mathcal{A}' = \mathcal{I} = \bigcup$
2. while $\min(\Delta\phi_{\mathcal{A}'}) < 0$ do
3. $\mathcal{A} = \mathcal{A} \cup \mathcal{A}'(\Delta\phi_{\mathcal{A}'} < 0)$
4. $\mathcal{A}' = \mathcal{I} \setminus \mathcal{A}$
5. $\Delta\phi_{\mathcal{A}'}^{n+1} = -(\mathbf{K}_{\phi\phi}^{n+1})_{\mathcal{A}'\mathcal{A}'}^{-1}(\mathbf{f}_{\phi}^{n+1})_{\mathcal{A}'}$
6. $\Delta\phi_{\mathcal{A}} = 0$
7. end
8. $\phi^{n+1} = \phi^n + \Delta\phi^{n+1}$

### 4.5. Adaptive element deletion

In the Lagrangian coordinates frame, the grid deforms with the material, which may engender severe element distortion, especially in the large-strain fracture problem involving discontinuous surfaces. Despite phase field representation of the fracture erasing the explicit discontinuity, the mesh distortion issue is not settled, but even worse. As displayed in Fig. 3 (a), the crack opening zone that should have no material is now filled with damaged elements. To delineate the crack opening outline, these elements with degraded stiffness will undergo large deformations. Moreover, the most severe mesh deformation arises precisely in these damaged zones. A consequential idea is that if these damaged elements can be properly removed, the mesh distortion dilemma facing large strain fracture will be significantly alleviated. Following this route, an adaptive mesh deletion technology for phase field modeling of large-strain fracture is employed. We first formulate a rule for element deletion, i.e., the elements that meet

$$\max\left(\phi_e^i\right) \geq \phi_c \tag{55}$$

will be removed from the solution domain (The brief algorithm procedures refer to [49]). From our experience, the setting $\phi_c \in [0.95, 0.98]$ is apposite. After eliminating the labeled elements, the crack opening profile is visualized, as illustrated in Fig. 3(b).

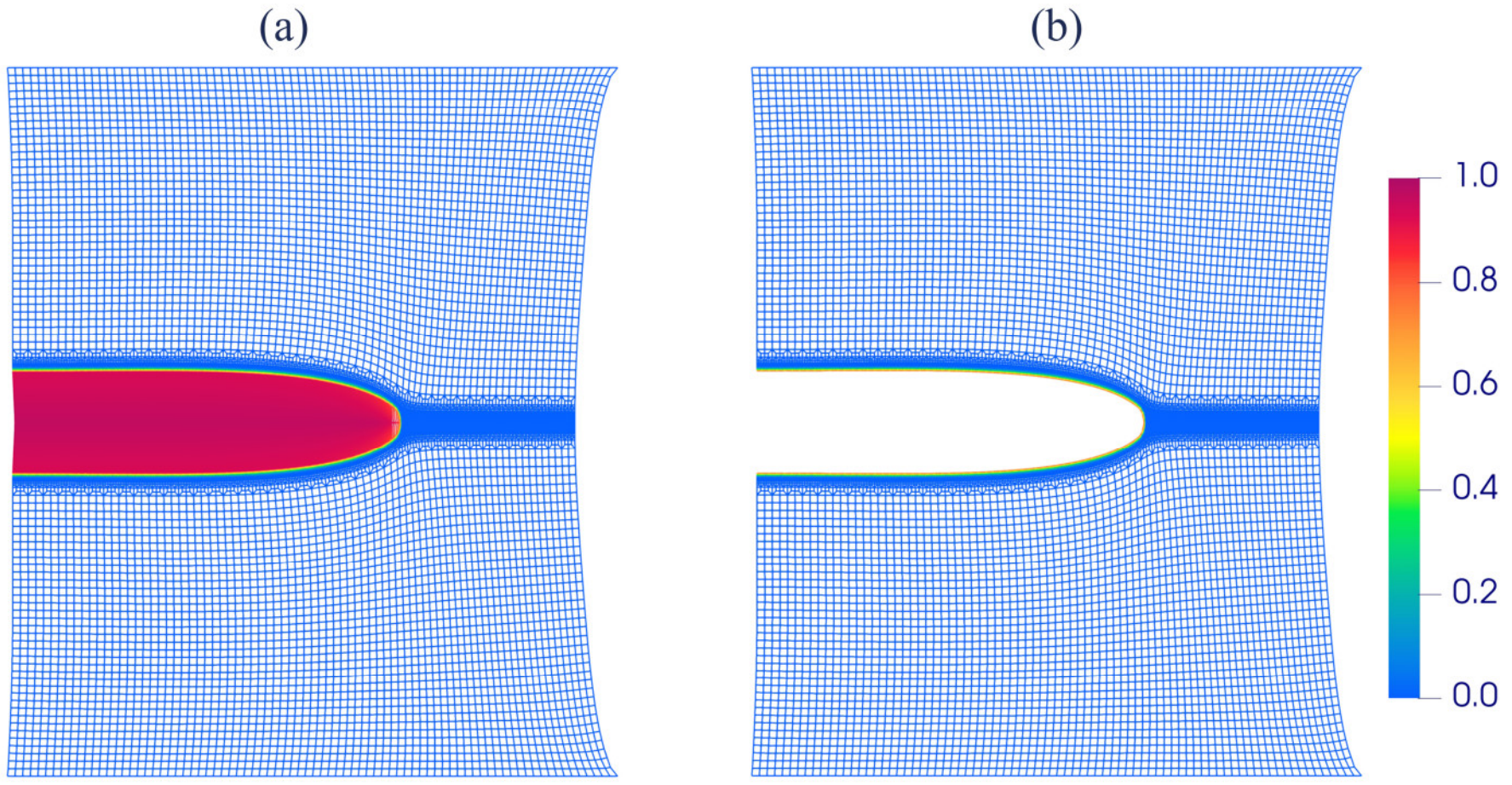

Fig. 3. The mesh configuration before (a) and after (b) the damaged element is deleted.

### 4.6. Solution schemes

Throughout the solving procedures, the pre-processing mesh generation module is developed on the foundation of open-source software packages termed ifem [59] and ameshref [60]. The proposed mixed displacement-pressure-phase field formulation for finite strain fracture of quasi-incompressible hyperelastic materials is realized by the parallel in-house MATLAB code and post-processed with ParaView [61]. To gain insight into the complete solution scheme, the basic procedures and core technologies of the proposed approach are recapitulated in Algorithm 2.

**Algorithm 2**

1. Generate FEM mesh using Q4(H8) elements.
2. Initialize $\phi_0$, $\boldsymbol{u}_0$, $p_0$ and $tol = 1\times10^{-4}$.
3. Create a pre-crack $\mathcal{C}$ by setting $\phi_{\mathcal{C}} = 1$.
4. For load step $n+1$, run

5. while ($Res^{u-p} \geq tol$ or $Res^{\phi} \geq tol$) do

Assemble $\boldsymbol{K}_{\boldsymbol{uu}}^{i+1}, \boldsymbol{K}_{\boldsymbol{up}}^{i+1}, \boldsymbol{K}_{pp}^{i+1}, \boldsymbol{f}_{\boldsymbol{u}}^{i+1}, \boldsymbol{f}_{\boldsymbol{u}}^{i+1}$ with frozen $\phi_{n+1}$.

Update $\boldsymbol{u}_{n+1}^{i+1}$, $p_{n+1}^{i+1}$ by solving

$$\begin{bmatrix} \boldsymbol{u}_{n+1}^{i+1} \\ p_{n+1}^{i+1} \end{bmatrix} = \begin{bmatrix} \boldsymbol{u}_{n+1}^{i} \\ p_{n+1}^{i} \end{bmatrix} + \begin{bmatrix} \boldsymbol{K}_{\boldsymbol{uu}}^{i+1} & \boldsymbol{K}_{\boldsymbol{up}}^{i+1} \\ \left[\boldsymbol{K}_{\boldsymbol{up}}^{i+1}\right]^T & \boldsymbol{K}_{pp}^{i+1} \end{bmatrix}^{-1} \begin{bmatrix} \boldsymbol{f}_{\boldsymbol{u}}^{i+1} \\ \boldsymbol{f}_{p}^{i+1} \end{bmatrix} \tag{56}$$

Assemble $\boldsymbol{K}_{\phi\phi}^{i+1}$, $\boldsymbol{f}_{\phi}^{i+1}$ with frozen $\boldsymbol{u}_{n+1}$, $p_{n+1}$.

Update $\phi_{n+1}^{i+1}$ according to Algorithm 1.

Assign

$$Res^{u-p} = \max\left\{\frac{\left\|\boldsymbol{f}_u^{i+1}\right\|}{\left\|\boldsymbol{f}_u^{i}\right\|}, \frac{\left\|\boldsymbol{f}_p^{i+1}\right\|}{\left\|\boldsymbol{f}_p^{i}\right\|}\right\}, \quad Res^{\phi} = \frac{\left\|\boldsymbol{f}_\phi^{i+1}\right\|}{\left\|\boldsymbol{f}_\phi^{i}\right\|} \tag{57}$$

with the Euclidean norm $\|\cdot\|$.

End do

6. Perform adaptive mesh deletion (optional).
7. Go into the next Load step and rerun step 5.

## 5. Validation via numerical examples

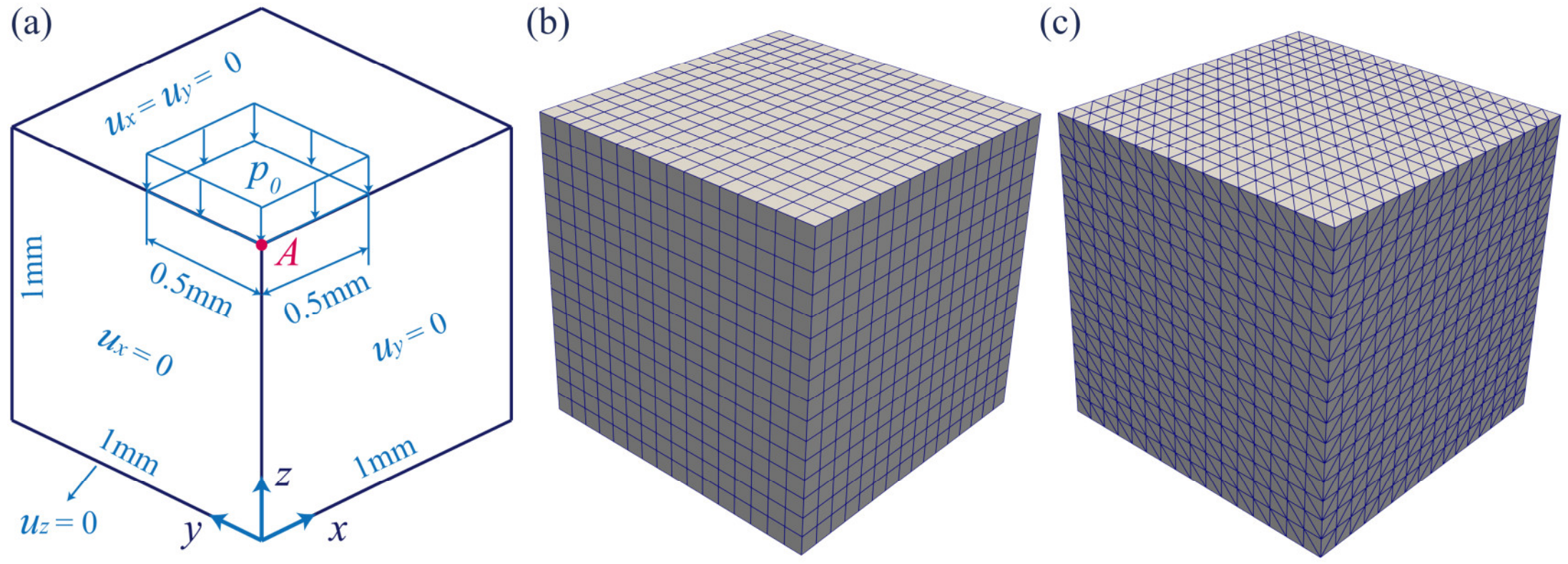

Fig. 4. 3D block under compression. (a) Quarter geometry and boundary conditions. The mesh configurations for Q1/P0 (b) and P2/P1 formulations (c).

In this section, the competence of the proposed framework for modeling the fracture of quasi-incompressible hyperelastic materials is validated by several numerical tests. The first example considers a well-known benchmark test of a 3-D block under compression to demonstrate the effectiveness of mixed Q1/P0 and P2/P1 formulations in coping with quasi-incompressible issues by quantitative comparison with the standard displacement solutions. The ensuing two examples focus on plane-strain problems, consisting of a bilateral notched specimen in tension rooted in the experiments of Hocine et al. and an amusing peeling test with weak interfaces. The last one is the more challenging tearing test of the 3-D sheet at large deformation, showing the superior performance of the proposed formulation and numerical treatment.

### 5.1. Quasi-incompressible block under compression

This example is a widely accepted benchmark test for evaluating the performance of numerical algorithms in approaching the incompressible limit. The strain energy function takes the same form as Reese et al.[62], viz. compressible Neo-Hookean model

$$\psi = \frac{\mu}{2}\left(tr\boldsymbol{C} - 3\right) - \mu \ln\left(J\right) + \frac{\kappa}{2}\left(\ln J^2\right) \tag{58}$$

Herein, the incompressibility is guaranteed by setting $\kappa \gg \mu$, such that the shear modulus $\mu = 80.194$ N/mm$^2$ and the bulk modulus $\kappa = 400889.806$ N/mm$^2$. Fig. 4 (a) presents the quarter geometric configuration and boundary conditions. A uniform pressure $p_0$ is imposed on a quarter square zone encompassing the top surface center (point $A$, see Fig. 4(a)). The meshes used for Q1/P0 ($16\times16\times16$) and P2/P1 ($[16\times16\times16]\times6$) formulations are demonstrated in Figs. 4(b) and (c), respectively. Besides, the identical mesh shown in Fig. 4(b) is also utilized to solve the standard displacement solution denoted by Q1S. By setting the pressure load $p_0 = 320$ N/mm$^2$ the displacement contours in deformed configuration for the above three schemes, Q1S, Q1/P0 and P2/P1 are presented together in Fig. 5. Thereinto, the standard Q1S exhibits over stiffening on account of the well-known volumetric locking, while Q1/P0 averts this issue, closely resembling the response of P2/P1 formulation in vision.

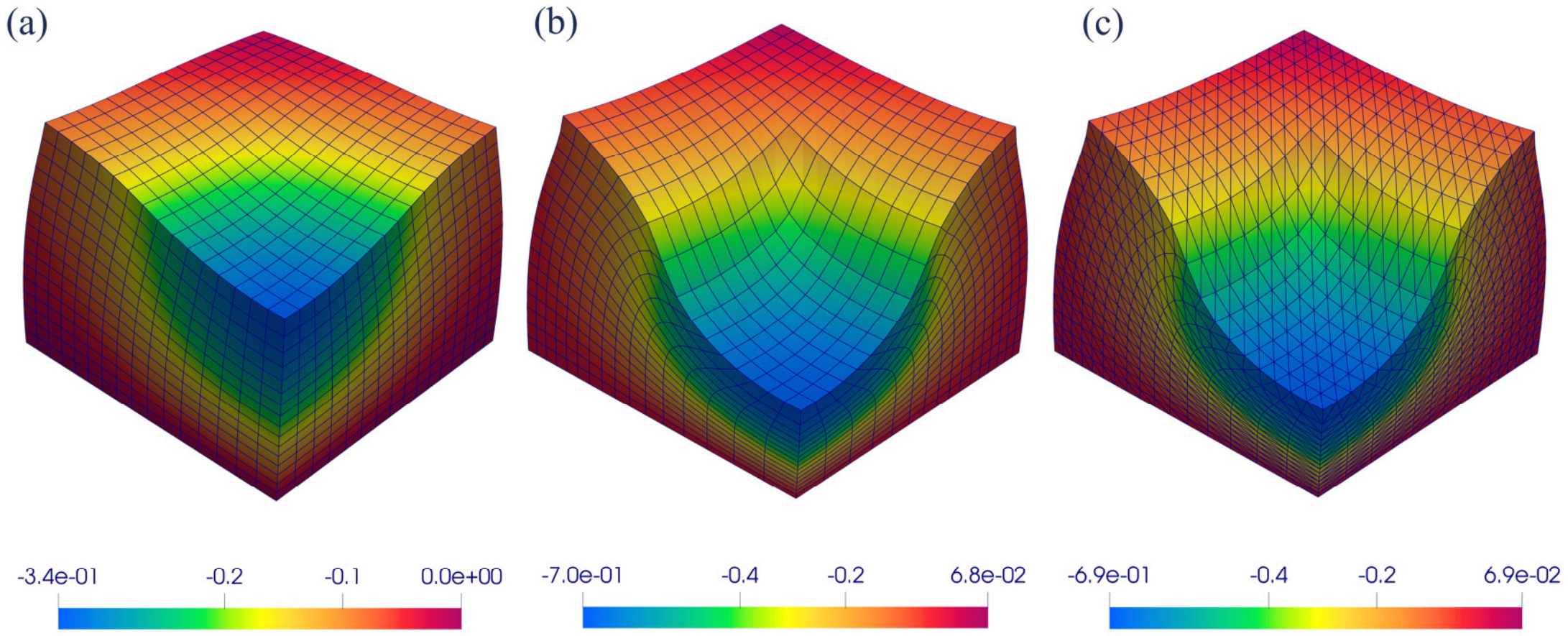

**Fig. 5.** Displacement contour plots of the 3D block under compression in deformed configuration for Q1S solution (a), mixed Q1/P0 formulation (b) and P2/P1 formulation (c).

For quantitative comparison, the relative z-displacement of point $A$ of the three schemes, known as the so-called compression levels, are exported. Fig. 6 (a) depicts the convergence of compression level with the number of elements per edge. As shown, the compression levels obtained by Q1/P0 and P2/P1 formulations are almost indistinguishable after the number of elements per edge exceeds 12. The result calculated by the Q1S scheme, by contrast, is significantly lower, even if the mesh is refined. A similar phenomenon is also reflected in Fig. 6 (b), which depicts the evolution of the compression level with loading for a mesh configuration of 16 elements per edge. We hereby underline that this benchmark test corroborates the efficacy of the basic framework for incompressible problems, thus laying a solid foundation for subsequent incorporating phase field.

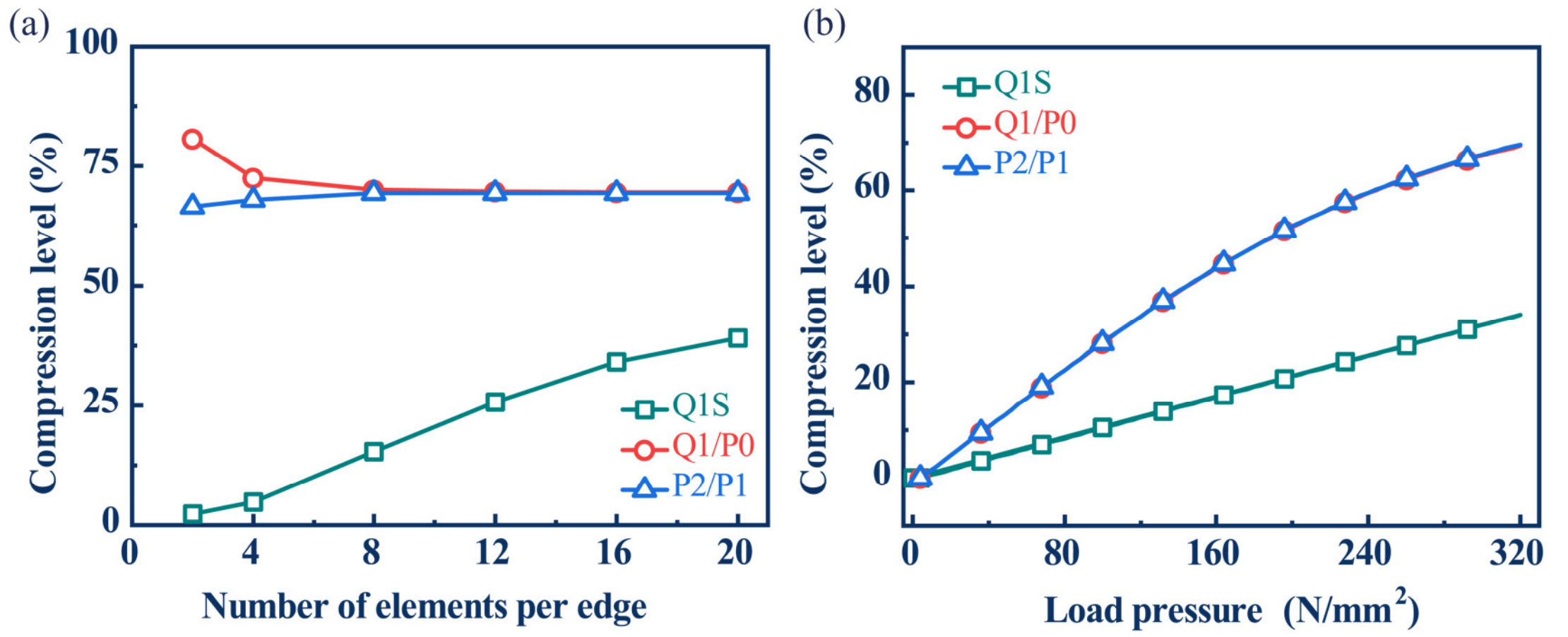

**Fig. 6.** Comparison of calculation results in three formulations. (a) Study on the convergence of compression level (point A) with the number of elements per edge. (b) Comparison for the evolution of compression level with the loading.

## 5.2. Double edge notched specimen in tension

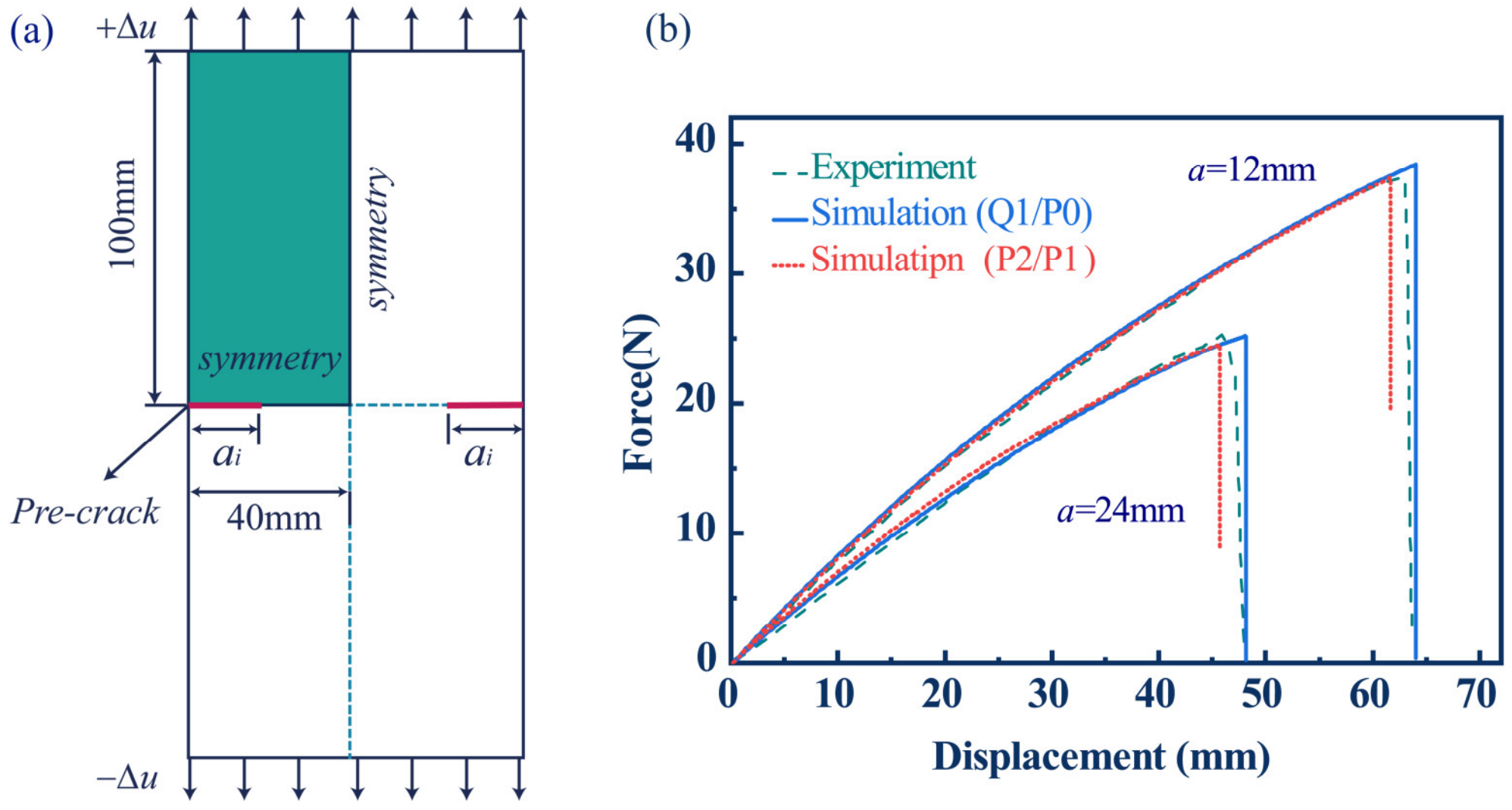

**Fig. 7.** Double edge notched specimen in tension. (a) Geometry and boundary conditions allowing for symmetry. (b) Comparison of simulated (Q1/P0 and P2/P1) and experimental force-displacement responses.

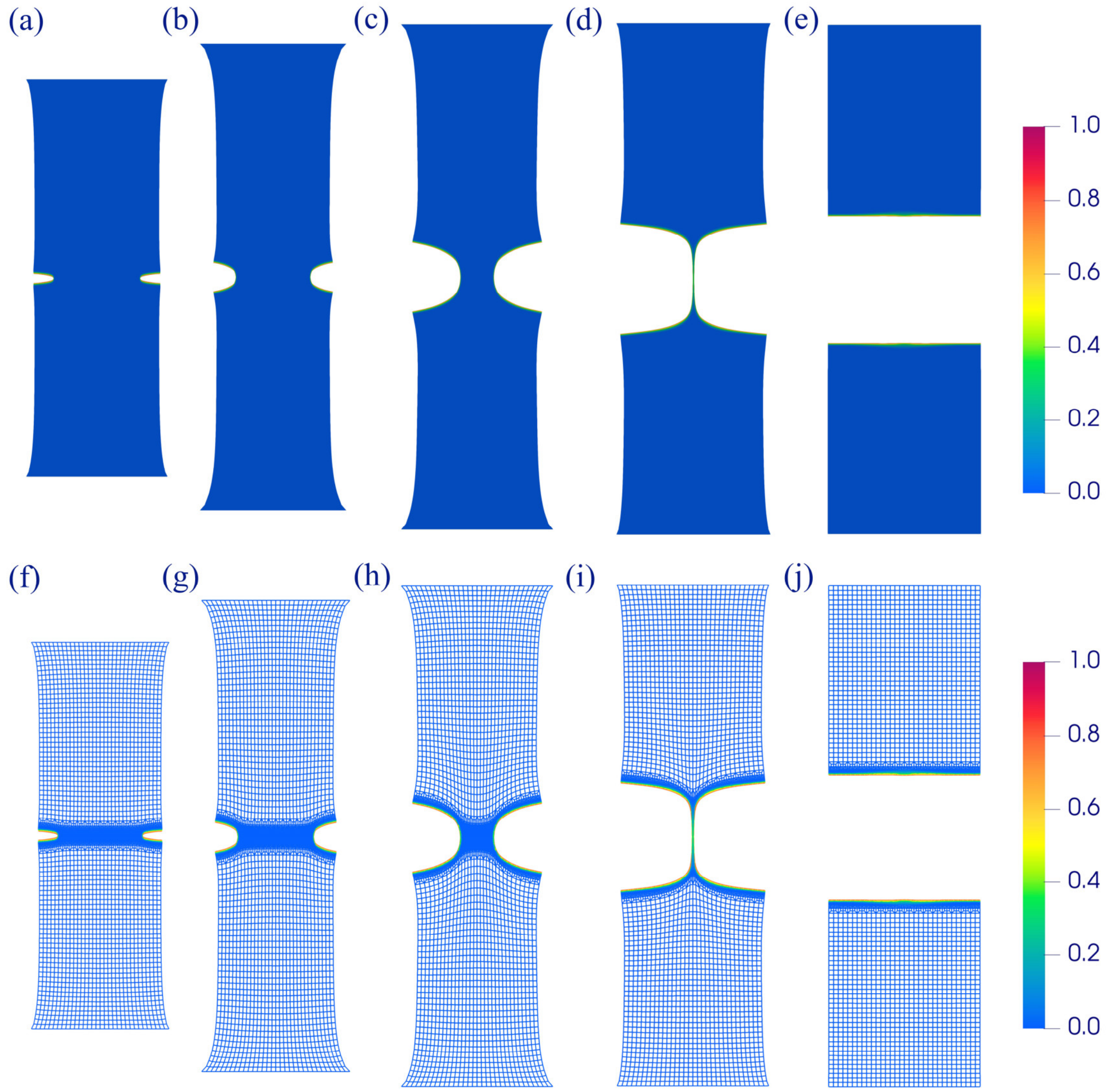

**Fig. 8.** Crack patterns and the corresponding meshes of the double-edge tensile specimen with the notch length of $a_i = 12$ mm in the deformed configuration. The snapshots (a)-(e) (or (f)-(j)) are at loading displacements of $\Delta u = 10$, $\Delta u = 40$, $\Delta u = 64.0263$, $\Delta u = 64.0295$ $\Delta u = 64.0296$ [mm].

This example aims to test the performance of the proposed mixed formulation for fracture modeling of incompressible hyperelastic materials by comparing with experiments [63]. Fig. 7 (a) presents the geometry and boundary conditions considering symmetry, aligning with the experiment of Hocine et al. The pre-crack length $a_i$ take 12 mm and 24mm, respectively. Using the Neo-Hookean model to characterize the mechanical response, the material parameters are set as: $\mu = 0.178$ N/mm$^2$, $v = 0.499$ and $G_c = 1.67$ N/mm. Following the previous reports, the phase-field regularization parameter is set to $l_0 = 1$ mm, the effective element size $h_f = 0.2$ mm and viscosity coefficient $\eta = 1\times10^{-3}$. With the above parameters, this example was simulated in terms of Q1/P0 and P2/P1 schemes, respectively, and the resulting force-displacement responses are plotted in Fig 7(b) together with the experimental ones. While the curves obtained by Q1/P0 and P2/P1 are in good agreement with the experimental measurements [63], the P2/P1 formulation breaks down when approaching the complete fracture state for the identical parameter settings.

To better visualize the crack opening morphology, the level set $\phi > 0.8$ is removed in ParaView[†]. The post-processed snapshots of the crack growth patterns and the corresponding Q1/P0 meshes at various fracture states are exhibited in Fig. 8 (see the supplemental material named mov1 for the complete crack evolution video), bearing a strong resemblance to those in the literature. We also derived the Cauchy stress and pressure fields with smooth distribution in the Q1/P0 formulation, and the contour of the Cauchy stress component $\sigma_{22}$ is demonstrated in Fig. 9. It is conspicuous that the crack opening outlines solved by the proposed scheme are smoother than past reports on this issue, by degenerating the incompressibility of damaged materials. Furthermore, Fig. 10 reports the volume evolution after eliminating the damage phase throughout the fracture process, corroborating that incompressibility is ensured for undamaged materials in the mixed Q1/P0 formulation. As for the P2/P1 scheme, the resultant crack patterns are indistinguishable from those of Q1/P0 (see **Appendix B**), except for the poorer robustness. For this reason, the P2/P1 scheme is shelved in the upcoming tests.

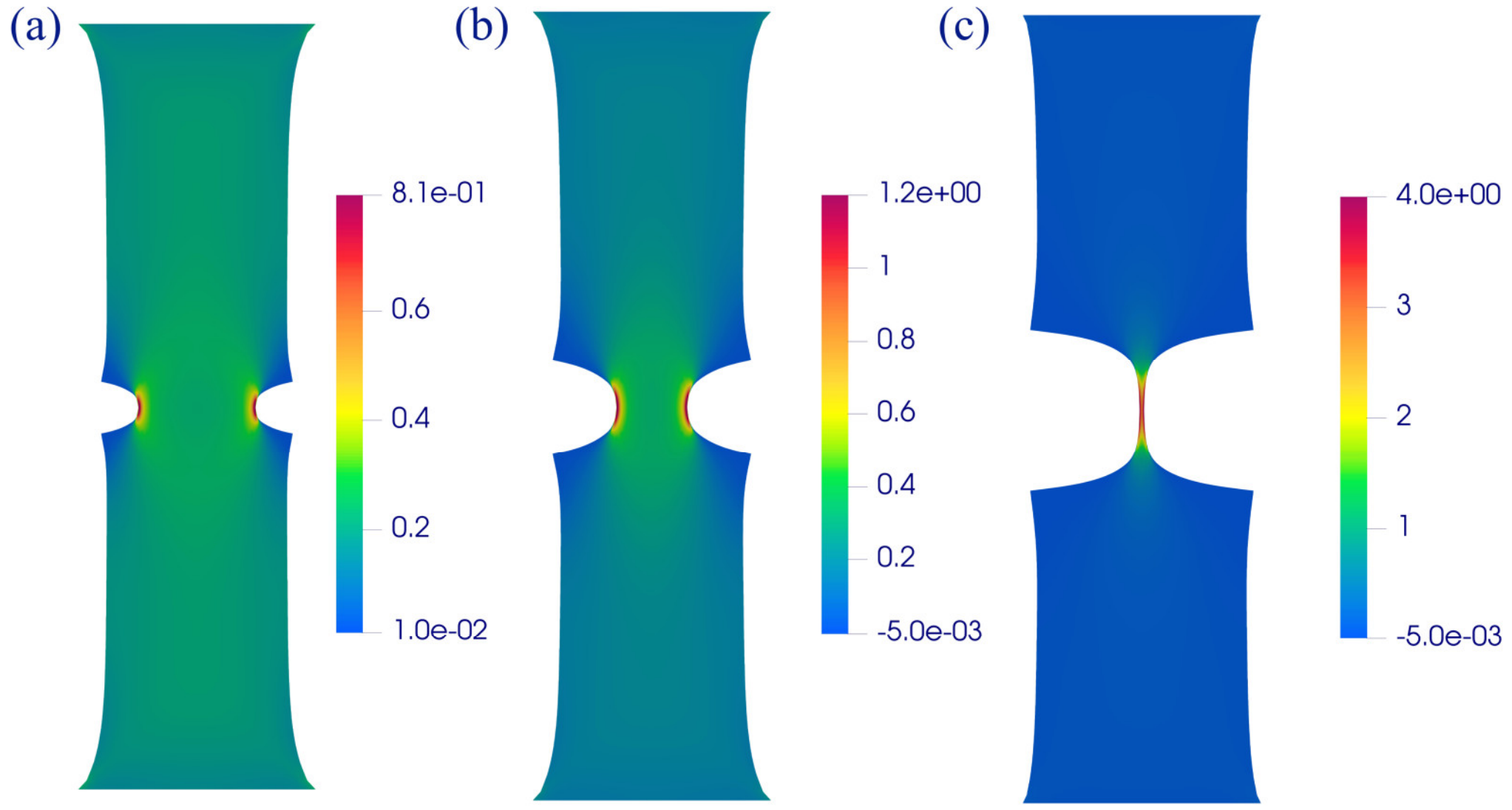

**Fig. 9.** Distributions of Cauchy stress $\sigma_{22}$ for Q1/P0 formulation at various loading states: $\Delta u = 40$, $\Delta u = 64.0241$, $\Delta u = 64.0293$ [mm].

[†] This operation is a post-processing technique, unlike the foregoing adaptive mesh deletion method.

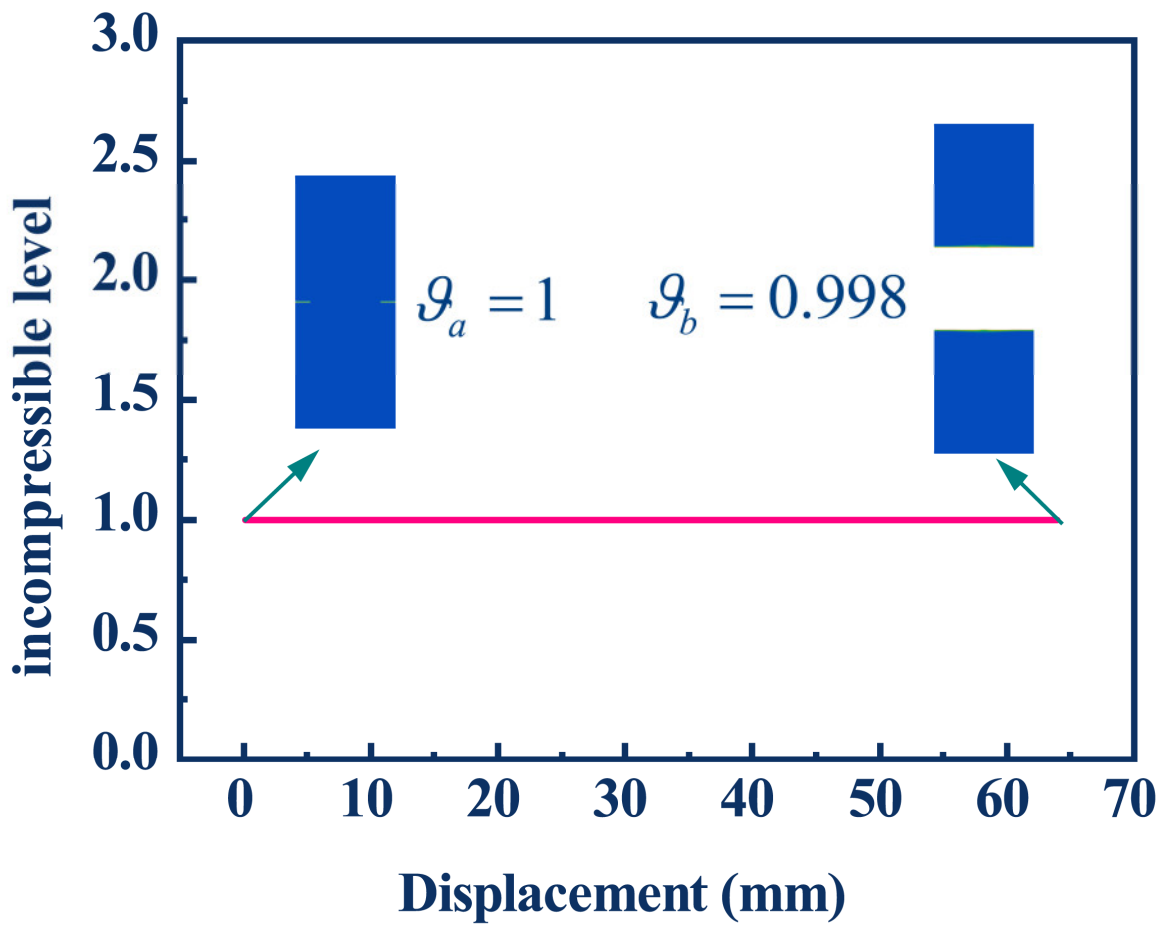

**Fig. 10.** Evolution of the incompressible level $\vartheta$ throughout the fracture process.

## 5.3. Two-dimensional peeling test

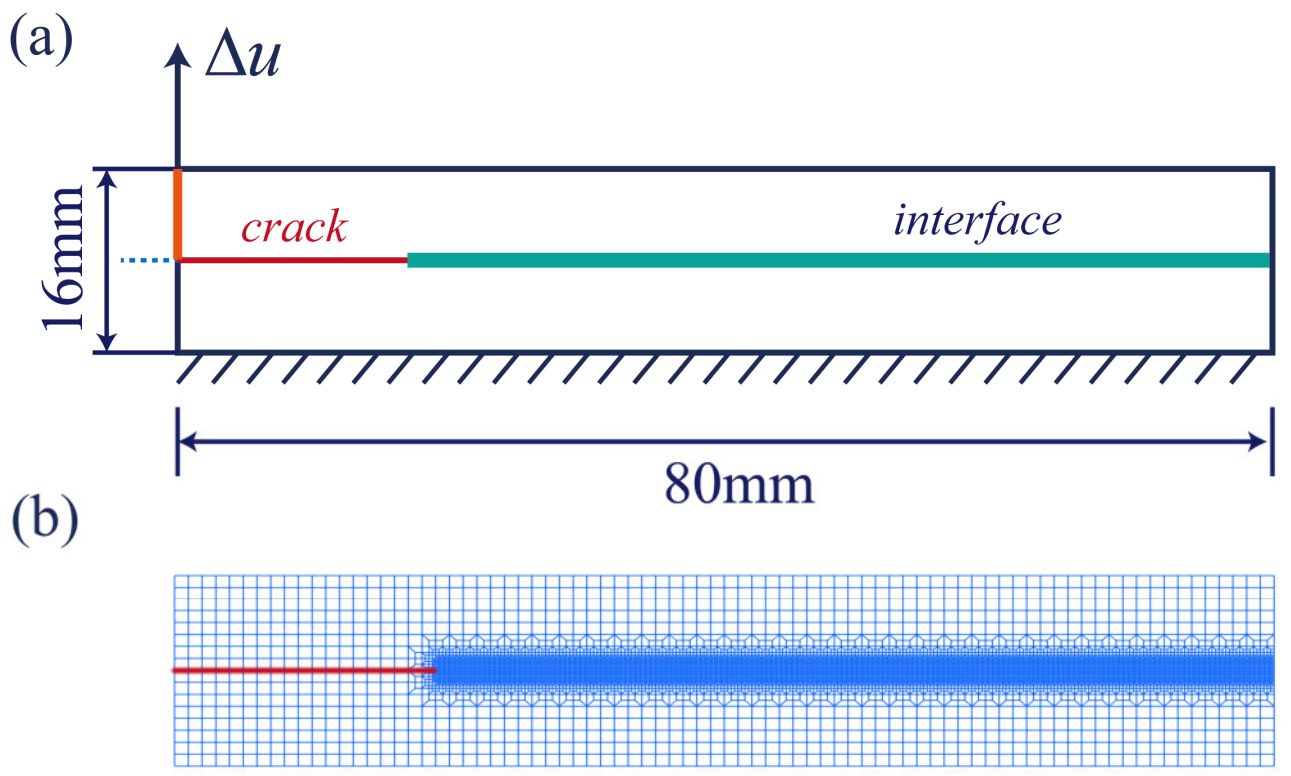

**Fig. 11.** Peeling test in soft materials containing the weak interface. (a) Geometry configuration and boundary conditions. (b) finite element meshes used for Q1/P0 formulation.

The design of this example is inspired by an intriguing gel peeling experiment [64], involving more complicated configurations and stronger nonlinearities. We consider a rectangular sheet of 80 mm length and 16 mm width, containing a weak interface with a width of 0.5mm at the middle height, as depicted in Fig.11 (a). The left edge is split into two segments by a notch of 20 mm in length at the interface. With the fixed bottom edge, the seed crack can be driven by imposing a vertical displacement loading $\Delta u$ on the upper segment of the left edge. Fig. 11 (b) presents the finite element meshes used for this example, where the elements in the vicinity of the interface are refined (effective element size of $h_f = l_0 / 5$). For the sake of testing the serviceability of the mixed formulation, we switch to a more realistic Mooney-Rivlin model (see Eq. 7 for the strain energy function), in which the material parameters are set as: $c_1 = 0.078$ N/mm$^2$ ,

$c_2 = 0.054$ N/mm$^2$, $v = 0.499$, $G_c^{bulk} = 1.78$ N/mm and $G_c^{interface} = G_c^{bulk} / 20$. The other requisite parameters are $l_0 = 0.5$ mm and $\eta = 1 \times 10^{-3}$.

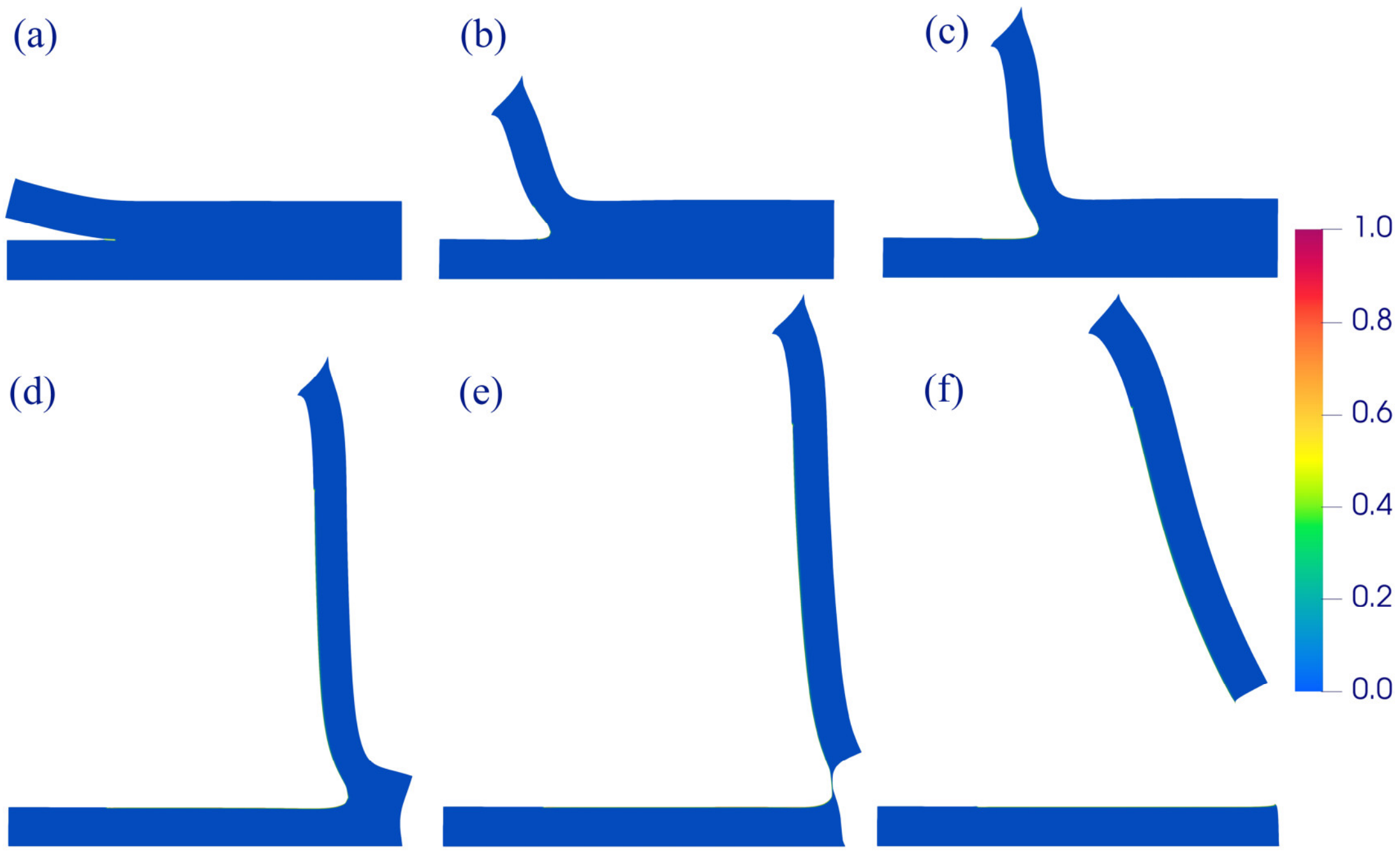

**Fig. 12.** Deformation states for the 2-D peeling test at loading displacements of $\Delta u = 4$ (a), $\Delta u = 27.62$ (b), $\Delta u = 36.24$ (c), $\Delta u = 86.06$ (d), $\Delta u = 94.1835$ (e), $\Delta u = 94.2633$ [mm] (f).

Fig. 12 depicts the entire peeling process from the pre-crack opening (Fig. 12(a)), germination (Fig. 12(b)), growth (Figs. 12(c)-e), and final separation at the weak interface (Fig. 12(f)). The evolution course of the crack patterns in the deformed configuration is quite similar to experimental observations. Then, we also extracted the pressure and Cauchy stress fields, and a snapshot of their smooth distributions at the loading of $\Delta u = 74.82$ [mm], is presented in Fig. 13. For better visualization, a complete landscape of the entire fracture process is placed in supplemental materials and termed mov2. Now, despite this example demonstrating the outstanding performance of the proposed mixed framework, more challenging three-dimensional fracture scenarios are still lacking, which yields motivation for the following test.

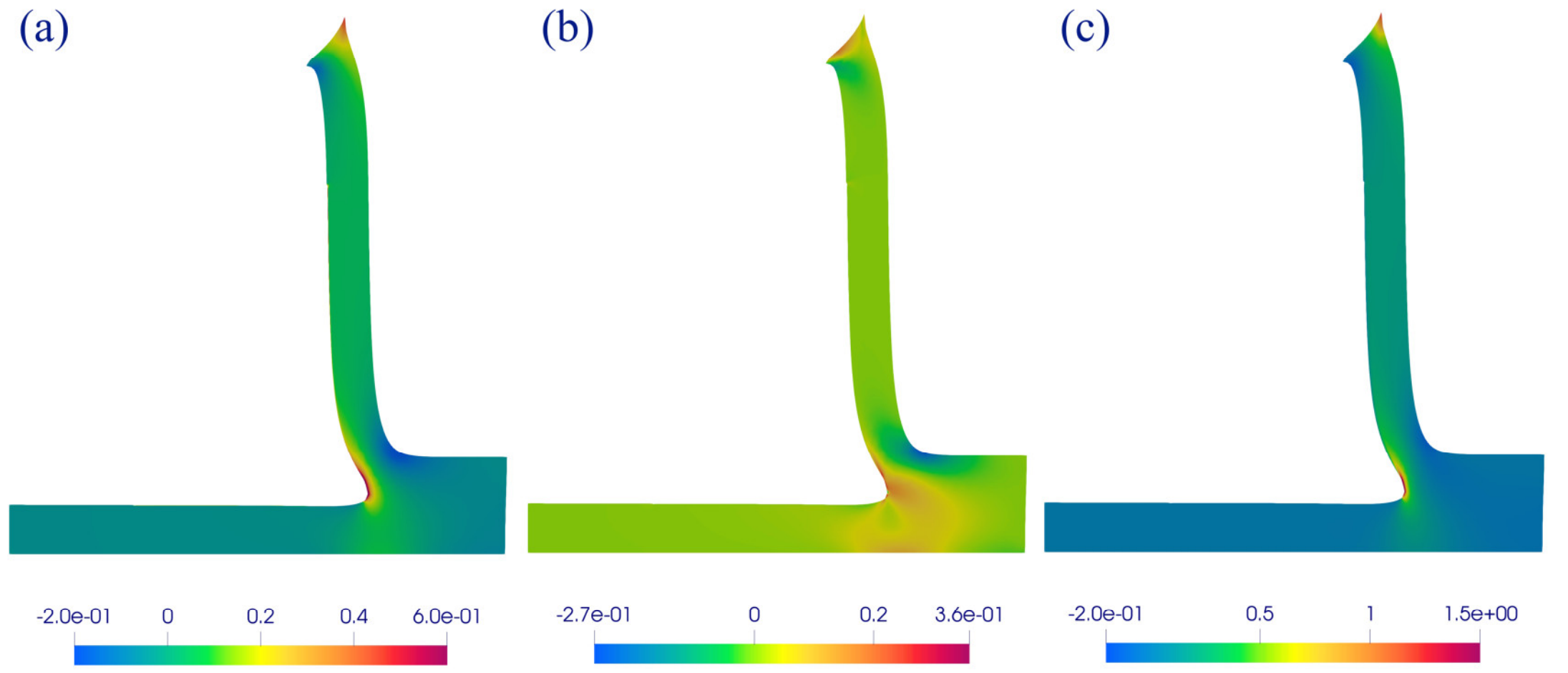

**Fig. 13.** Distributions of pressure $p$, Cauchy stress components $\sigma_{11}$ and $\sigma_{22}$ for Q1/P0 formulation at the loading of $\Delta u = 74.82$ [mm].

## 5.4. Tearing test in three-dimensional sheet

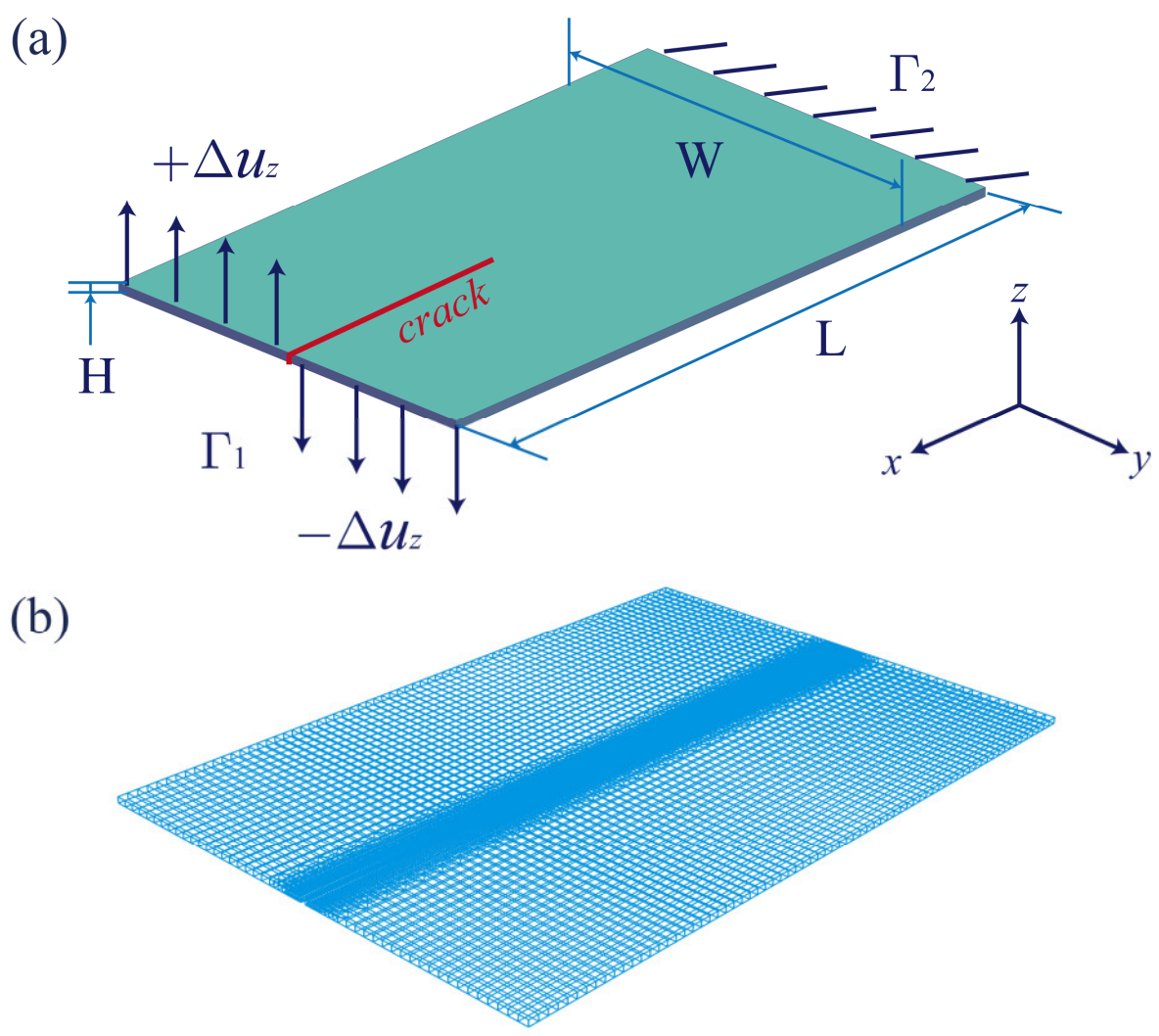

**Fig. 14.** Tearing test in 3-D sheet. (a) Geometrical setup and boundary conditions. (b) finite element mesh used for Q1/P0 formulation.

The idea of this example stems from the classical out-of-plane tearing experiment [65], a bit like the test of Yin et. al [66]. We consider a 3-D rubber sheet of 48 mm length ($L$), 36 mm width ($W$) and 1 mm thickness ($H$), containing a seek crack of 12 mm length, as illustrated in Fig. 14(a). Constraining the right side $\Gamma_2$, the sheet is torn apart by exerting antisymmetric z-displacements to the two left sides split by the crack. The finite element mesh used for this test is presented in Fig. 14(b), where the a priori crack path is refined to meet the effective element size of 0.125mm. We assume that the strain energy function is portrayed by the Neo-Hookean model for simplicity. The required calculating parameters as $\mu = 0.203$ N/mm$^2$, $G_c = 0.72$ N/mm, $\nu = 0.499$, $l_0 = 0.4$ mm and $\eta = 1\times10^{-3}$.

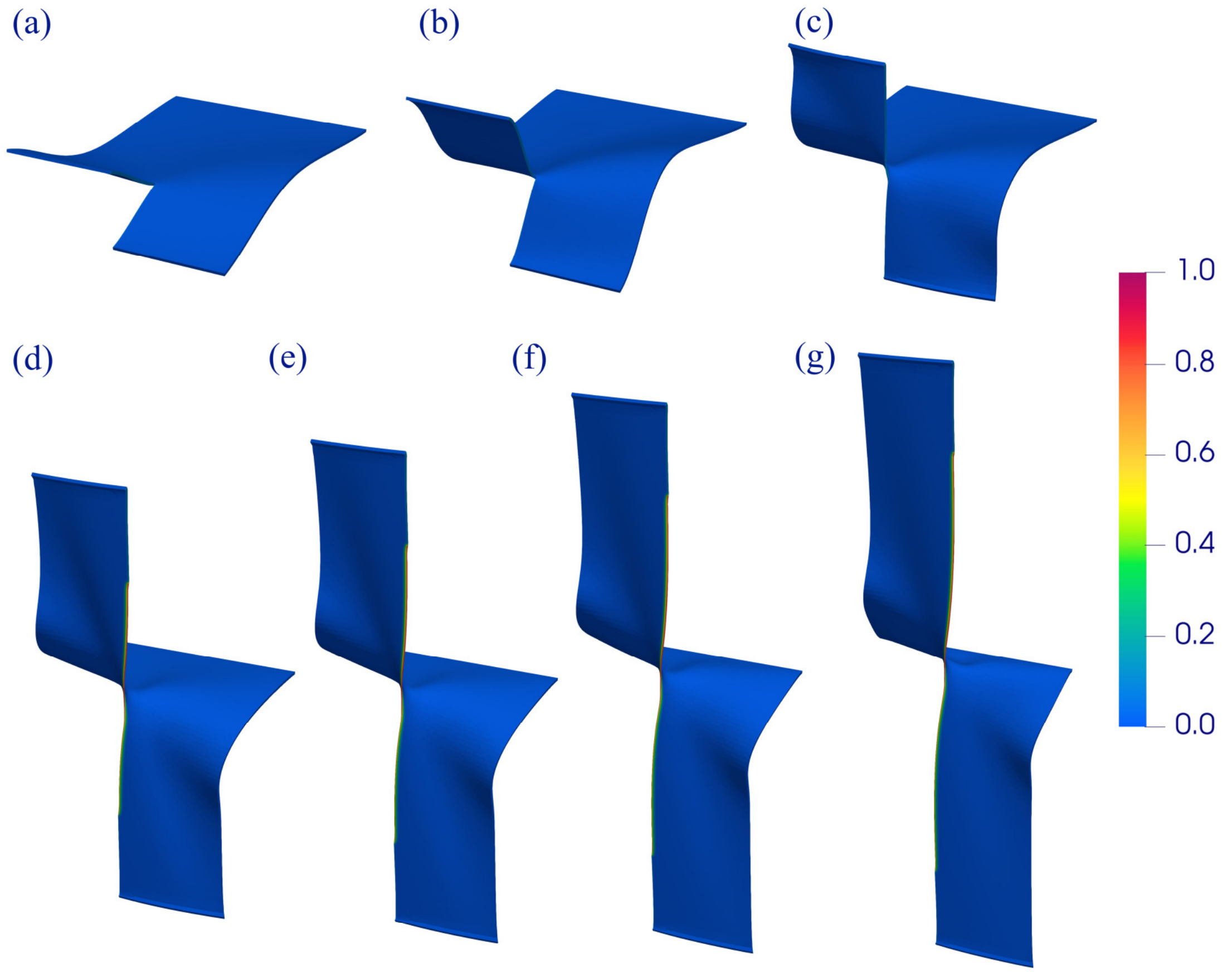

**Fig. 15.** Snapshots of 3-D tear test in deformed configuration at various loading states: $\Delta u = 9.68$ (a), $\Delta u = 16.0$ (b), $\Delta u = 30.08$ (c), $\Delta u = 64.08$ (d), $\Delta u = 80.08$ (e), $\Delta u = 87.68$ (f), $\Delta u = 100.04$ [mm] (g).

Fig. 15 displays the deformation history of out-of-plane tearing through a series of snapshots in various loading states. The sheet first undergoes a purely elastic deformation phase in the early stage of loading, see Figs. 15(a) and (b). As the loading approaches $\Delta u \sim 30.08$ mm, the damage initiates at the forefront of the crack opening outline (Fig. 15 (c)). After that, the loading is further increased to ignite the tearing process, as depicted in Figs 15 (d)-(g). We also present the distribution of other physical fields in an arbitrarily selected loading state, comprising pressure $p$ and the magnitudes of displacement $\|\boldsymbol{u}\|$ and Cauchy stress $\|\boldsymbol{\sigma}\|$, as shown in Fig. 16. A complete tearing animation termed mov3 is provided in supplemental materials, close resemblance to the experimental scene [66]. As a final note, we hereby illuminate that adaptive mesh deletion technology is vital to realizing this tearing test, although it is optional in the previous examples.

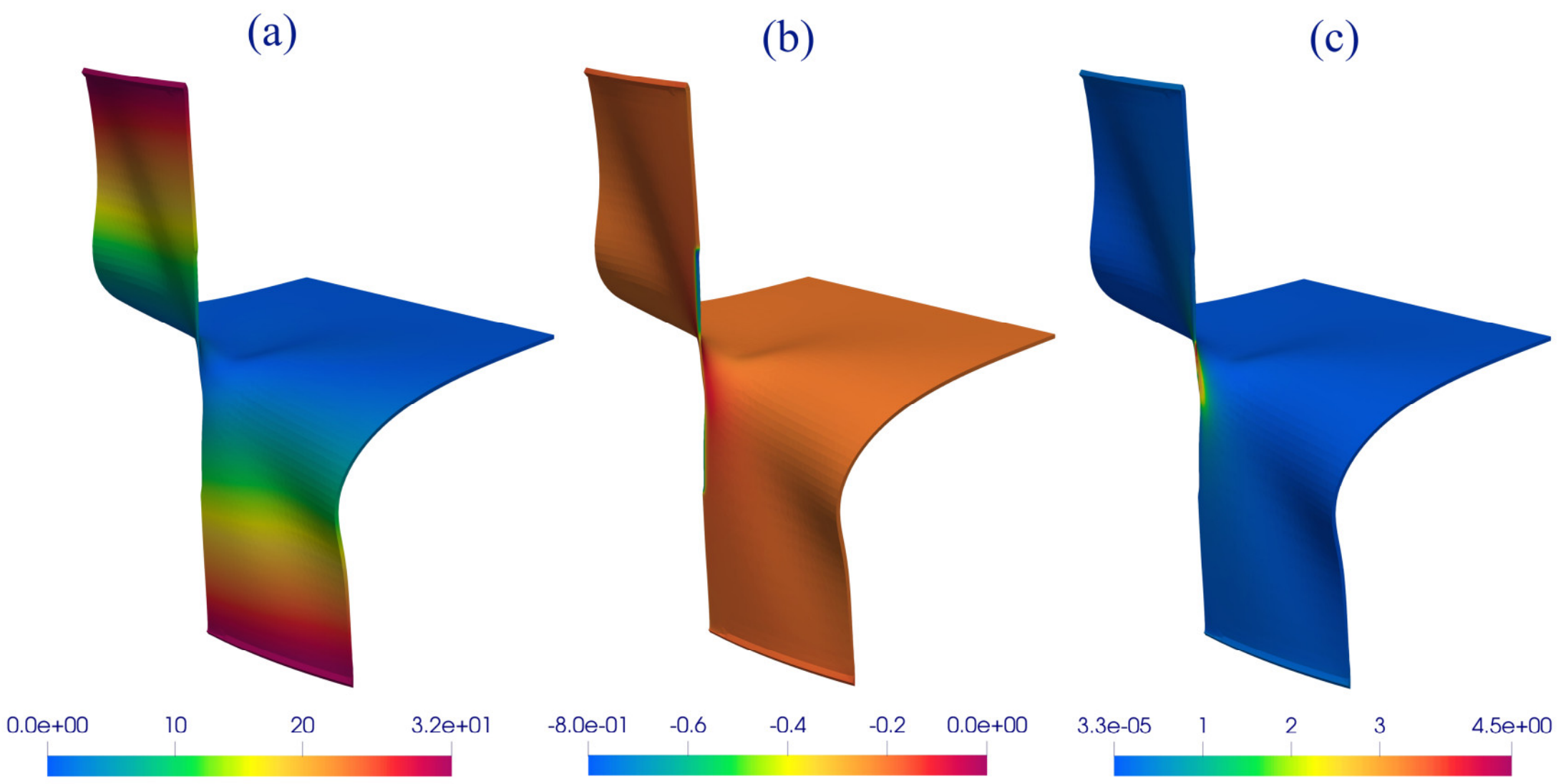

**Fig. 16.** Tear test in 3-D sheet. Distributions of the magnitude of displacement $\|\boldsymbol{u}\|$, pressure $p$ and the magnitude of Cauchy stress $\|\boldsymbol{\sigma}\|$ at the loading of $\Delta u = 58$ [mm].

## 6. Conclusion

In the finite strain regime, we have posed a novel mixed displacement-pressure-phase field framework for addressing the fracture issues of nearly incompressible hyperelastic materials. Based on the modified perturbed Lagrangian approach, a multi-field variational form is derived. The concrete numerical implementations of the proposed formulation are also presented, and its excellent performance is demonstrated via four numerical examples. For conciseness, we summarize the main contributions of this study as follows

(i) The innate contradiction between incompressibility and the diffuse crack opening at finite deformation was unveiled.

(ii) A scheme that loosens incompressible constraint of the damaged phase without affecting the intact material was proposed to resolve the underlying conflict.

(iii) In the numerical aspect, an adaptive mesh deletion technology is developed for the current issue, alleviating the mesh distortion under large deformations, especially in 3-D scenarios.

(iv) In terms of precision and robustness, the classical Q1/P0 scheme performs better in the fracture modeling of nearly incompressible hyperelastic materials than the higher-order P2/P1 formulation.

To sum up, the proposed mixed framework eliminates the barriers set by incompressibility for phase field crack opening, thus initiating a new avenue for finite strain fracture in nearly incompressible cases. The current formulation will be further extended to study dynamic fracture issues with richer crack morphologies in future work.

## ACKNOWLEDGMENTS

We acknowledge support from the National Natural Science Foundation of China (Grant No. 51890872 and 51790500).

## Appendix A. Incompressibility frustrates crack opening

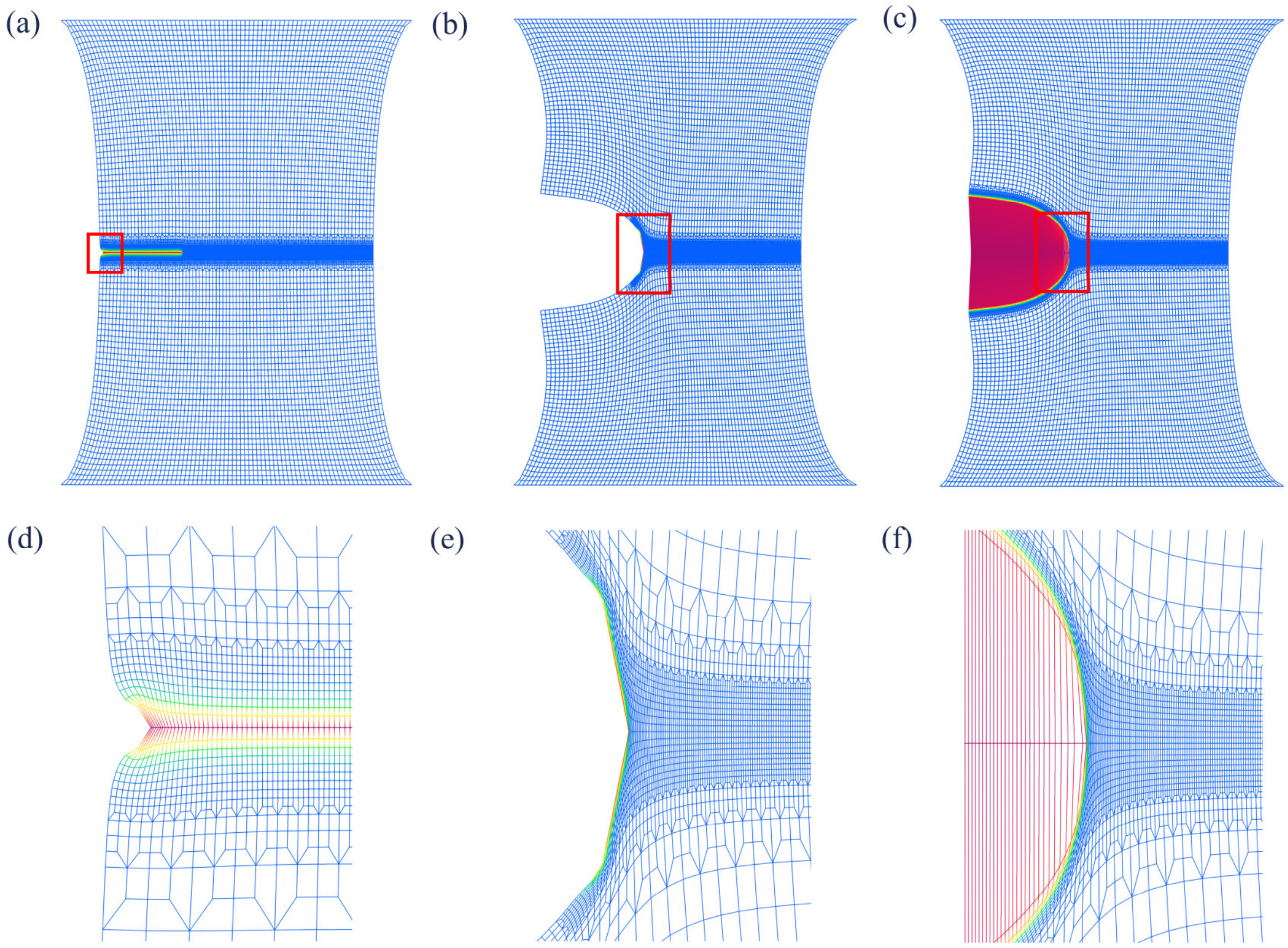

**Fig. A1.** Incompressible damaged phase obstructing cracks opening with diffuse (a) and geometric (b) pre-crack, respectively. (c) The results obtained by the proposed framework. (d)-(e) are the partial enlargements corresponding to (a)-(c).

Based on the geometric configuration illustrated in Fig. 1, the consequences of imposing incompressibility on damaged materials are demonstrated here (material parameters refer to section 5.2). We first consider diffuse pre-crack, i.e., set $\phi = 1$. As shown in Fig. A1(a) and its partial enlarged Fig. 1(d), the pre-crack cannot open due to the limited volume change of damaged elements. Then, Fig. A1(b) presents the results of using geometric pre-crack. Although the crack opening was unhindered, the meshes at the crack tip were severely deformed, causing the crack outline to be unsmooth (see Fig. A1(e)). Subsequent calculations even encountered mesh distortion issues, resulting in numerical collapse. For comparison, we also demonstrate the result of relaxing the incompressible constraint of damaged materials, see Fig. A1(c). The crack opening profile is very similar to the physical one, but there is no mesh distortion at the crack tip, as depicted in Fig. A1(f). Given that the damaged material has almost no substantial contribution to the calculation system, it is sound to relax the incompressibility with the damage field.

## Appendix B. Recalculated test 1 with P2/P1 formulation

We recalculated test 1 using the high-order P2/P1 formulation. Fig. A2 exhibits the snapshots of the crack patterns and the corresponding P2/P1 meshes at various loading states. Evidently, the results obtained by P2/P1 are almost identical to that of the Q1/P0 scheme, which is also validated by the force-displacement curve in Fig. 7. From the numerical perspective, however, the P2/P1 scheme is less robust. With the same parameter settings, such a formulation cannot reach the final complete break state despite the higher computational cost.

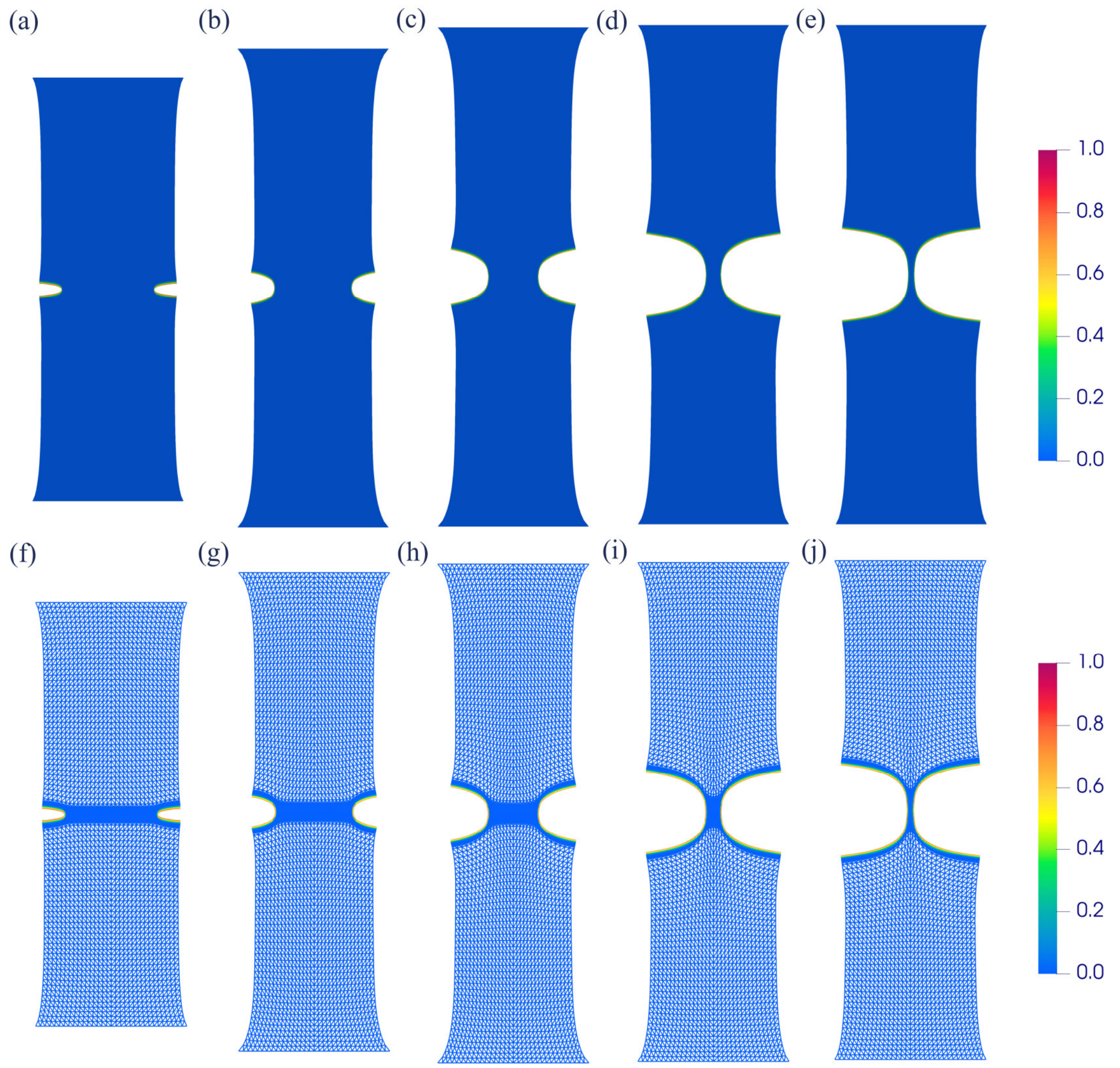

**Fig. A2.** Crack patterns and the corresponding meshes of test 1 recalculated by P2/P1 formulation. The snapshots (a)-(e) (or (f)-(j)) are at loading displacements of $\Delta u = 10$, $\Delta u = 40$, $\Delta u = 60.2573$, $\Delta u = 61.7294$ $\Delta u = 61.7296$ [mm].